\documentclass[12pt,journal,letterpaper,onecolumn]{IEEEtran2}
\usepackage{setspace}
\usepackage{graphicx}
\graphicspath{{pics/}{figs/}}
\usepackage{graphicx}
\usepackage{subfigure}
\usepackage{amssymb}
\usepackage{amsmath}
\def\C{\mathbb{C}}

\def\D{\mathbb{D}}
\newtheorem{defn}{\noindent $\mathbf{Definition}$}[section]
\newtheorem{lem}[defn]{$\mathbf{Lemma}$}

\newtheorem{thm}[defn]{$\mathbf{Theorem}$}
\newtheorem{cor}[defn]{$\mathbf{Corollary}$}
\def\thesection{\arabic{section}}
\def\thesubsection{\thesection.\arabic{subsection}}
\def\baselinestretch{1.3}
\begin{document}
%
\def\baselinestretch{1}
\title{Optimization of Surface Registrations using Beltrami Holomorphic Flow}
\author{$^\dag$L.M. Lui\thanks{$^\dag$ Lok Ming Lui, Department of Mathematics, Harvard University and UCLA\{malmlui@math.harvard.edu\}}, T.W. Wong, W. Zeng, X.F. Gu, P.M. Thompson, T.F. Chan, S.T. Yau}


\maketitle
\def\baselinestretch{1.3}
\begin{abstract}
In shape analysis, finding an optimal 1-1 correspondence between surfaces within a large class of admissible bijective mappings is of great importance. Such process is called surface registration. The difficulty lies in the fact that the space of all surface diffeomorphisms is a complicated functional space, making exhaustive search for the best mapping challenging. To tackle this problem, we propose a simple representation of bijective surface maps using Beltrami coefficients (BCs), which are complex-valued functions defined on surfaces with supreme norm less than 1. Fixing any 3 points on a pair of surfaces, there is a 1-1 correspondence between the set of surface diffeomorphisms between them and the set of BCs. Hence, every bijective surface map can be represented by a unique BC. Conversely, given a BC, we can reconstruct the unique surface map associated to it using the Beltrami Holomorphic flow (BHF) method. Using BCs to represent surface maps is advantageous because it is a much simpler functional space, which captures many essential features of a surface map. By adjusting BCs, we equivalently adjust surface diffeomorphisms to obtain the optimal map with desired properties. More specifically, BHF gives us the variation of the associated map under the variation of BC. Using this, a variational problem over the space of surface diffeomorphisms can be easily reformulated into a variational problem over the space of BCs. This makes the minimization procedure much easier. More importantly, the diffeomorphic property is always preserved. We test our method on synthetic examples and real medical applications. Experimental results demonstrate the effectiveness of our proposed algorithm for surface registration.
\end{abstract}


\thispagestyle{plain}

\section{Introduction}
Surface registration is a process of finding an optimal 1-1 correspondence between surfaces satisfying certain constraints. It is of great importance in different research areas, such as computer graphics and medical imaging. For example, in medical imaging, surface registration is always needed for statistical shape analysis, morphometry and the processing of signals on brain surfaces (e.g., denoising and filtering). In many cases, a surface must be non-rigidly aligned with another surface, while matching various features lying on both surfaces. Finding an optimal surface registration that best matches the required constraints is difficult, especially on convoluted surfaces such as the human brain. It is therefore necessary to develop an effective algorithm to compute the best surface registration.

In order to obtain the best 1-1 correspondence between two surfaces, an optimized surface registrations is often required. {\it Optimization of surface registrations} is the process of selecting an optimal surface diffeomorphism within a large class of admissible smooth mappings to best satisfy certain properties. It can usually be formulated as a variational problem in the form:
\begin{equation}
\min _{f \in \mathbb{F}_{\mathrm{Diff}}} E_0(f)
\end{equation}
\noindent where $\mathbb{F}_{\mathrm{Diff}} = \{f\colon S_1\to S_2\colon \mbox{$f$ is a diffeomorphism}\}$ is the space of all diffeomorphisms from surface $S_1$ to surface $S_2$.

Solving this type of variational problem is generally difficult, since the space of all surface diffeomorphisms $\mathbb{F}_{\mathrm{Diff}}$ is a complicated functional space. For instance, $\mathbb{F}_{\mathrm{Diff}}$ is inherently infinite dimensional and has no natural linear structure. Constructing an efficient optimization scheme in such space that guarantees to obtain a minimizer is a big challenge, and a loss of bijectivity of the surface maps (overlapping) is often observed during the optimization process. To solve this problem, it is necessary to develop a simple representation of surface diffeomorphisms which helps to simplify the optimization procedure.

In this paper, we propose a simple representation of surface diffeomorphisms using Beltrami coefficients (BCs). The BCs are any complex-valued functions defined on surfaces with $L^\infty$-norm strictly less than 1. Fixing any 3 points on a pair of surfaces, there is a one-to-one correspondence between the set of surface diffeomorphisms and the set of BCs. Hence, every bijective surface map can be represented by a unique BC. Conversely, given a BC, we propose to reconstruct the unique surface map associated to it using the Beltrami Holomorphic flow (BHF) method introduced in this paper. The BHF formulates the variation of the surface maps under the variation of BCs. Hence, variational problems of surface diffeomorphisms can be easily reformulated into variational problems of BCs in the form:
\begin{equation}
\min _{\mu \in \mathbb{F}_{\mathrm{BC}}} E(\mu)
\end{equation}
\noindent where $\mathbb{F}_{\mathrm{BC}} = \{\mu\colon S_1\to \mathbb{D} \colon ||\mu||_{\infty} <1 \}$ is the set of BCs.

The space of BCs is a much simpler functional space that captures the essential features of surface maps. There are no restrictions that BCs have to be 1-1, surjective or satisfy some constraints in their Jacobians. By adjusting BCs, we can adjust surface registrations accordingly using BHF to obtain surface maps with the desired properties. This greatly simplifies the minimization procedure. More importantly, the surface maps obtained is guaranteed to be diffeomorphic (bijective and smooth) during the optimization process. We have applied our proposed algorithm on synthetic examples and real medical applications for surface registration, which demonstrate the effectiveness of our proposed method.

In summary, our work contributes to the following three aspects:
\begin{itemize}
\item We propose a simple representation of surface diffeomorphisms to facilitate the optimization of surface registrations.
\item We develop a reconstruction algorithm of the surface diffeomorphism from a given BC, using BHF. This completes the representation scheme and allows us to move back and forth between BCs and surface diffeomorphisms.
\item  With BHF, we formulate variational problems of surface maps into variational problems of BCs. This greatly simplifies the optimization procedure.
\end{itemize}

\begin{figure*}[ht]
\centering
\includegraphics[height=2in]{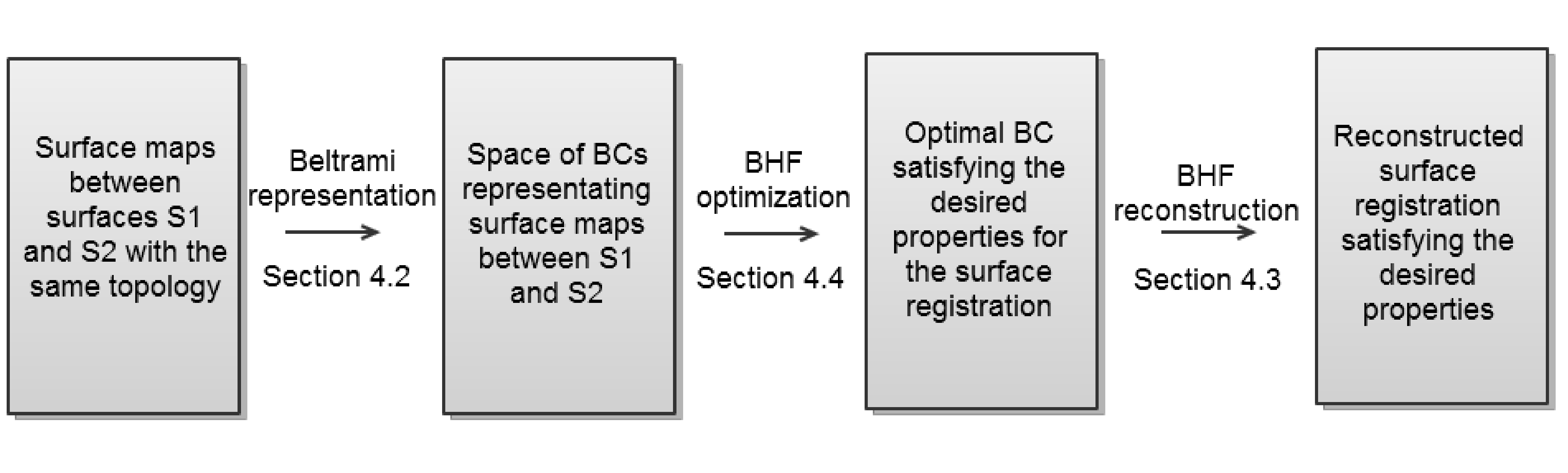}
\caption{A flow chart summarizing the framework proposed in this paper. \label{fig:BHFflowchart}}
\end{figure*}

A flow chart summarizing the framework proposed in this paper is shown in Figure \ref{fig:BHFflowchart}.

\section{Previous Work}
Surface registration has been extensively studied by different research groups. Most methods compute the optimal surface registration by minimizing certain kinds of energy functionals. In this section, we briefly describe some related methods commonly used.

Conformal surface registration has been widely studied to obtain smooth 1-1 correspondence between surfaces and minimize angular distortions \cite{Angenent}\cite{Gu1}\cite{Haker2}\cite{Hurdal2}\cite{Hurdal}\cite{Ju}\cite{Levy}\cite{Gu2}.. Conformal maps are usually computed using variational approaches to minimize some energy functionals, such as the harmonic energy \cite{Gu1} and the least square energy based on the Cauchy-Riemann equation \cite{Levy}. A 1-1 correspondence between surfaces can be obtained in the optimal state. However, the above registration cannot map anatomical features, such as sulcal landmarks, consistently from subject to subject.

To obtain a surface registration that matches important landmark features, landmark-based diffeomorphisms are often used. Optimization of surface diffeomorphisms by landmark matching has been extensively studied. Gu et al. \cite{Gu1} improves a conformal parameterization by composing an optimal M\"{o}bius transformation so that it minimizes a landmark mismatch energy. The resulting parameterization remains conformal, although features cannot be perfectly matched.  Wang et al. \cite{Lui07}\cite{Wang05} proposed a variational framework to compute an optimized conformal registration which aligns landmarks as well as possible. However, landmarks are not matched exactly and diffeomorphisms cannot be guaranteed when there is a a large amount of landmark features. Durrleman et al. \cite{Durrleman1}\cite{Durrleman2} developed a framework using currents, a concept from differential geometry, to match landmarks within surfaces across subjects, for the purpose of inferring the variability of brain structure in an image database. Landmark curves are not perfectly matched. Tosun et al. \cite{Tosun} proposed a more automated mapping technique that attempts to align cortical sulci across subjects by combining parametric relaxation, iterative closest point registration, and inverse stereographic projection. Glaun\`{e}s et. al \cite{glaunes}\cite{JoshiMiller} proposed to generate large deformation diffeomorphisms of a sphere onto itself, given the displacements of a finite set of template landmarks. The diffeomorphism obtained can better match landmark features. Lui et al. \cite{Lui10} proposed to compute shape-based landmark matching registrations between brain surfaces using the integral flow method. The one parameter subgroup within the set of all diffeomorphisms is considered and represented by smooth vector fields. Landmarks can be perfectly matched and the correspondence between landmark curves is based on shape information. Leow et al. \cite{leow:neuroimage2005} proposed a level-set-based approach for matching different types of features, including points, 2D and 3D curves represented as implicit functions. These matching features in the parameter domain were then pulled back onto surfaces to compute correspondence fields. Later, Shi et al. \cite{Shi07_2} computed a direct harmonic mapping between two surfaces by embedding both surfaces as the level-set of an implicit function, and representing the mapping energy as a Dirichlet functional in 3D volume domains. Although such an approach can incorporate landmark constraints, it has not been proven to yield diffeomorphic mappings.

Since there may not be well-defined landmarks on surfaces, some authors proposed driving features into correspondence based on shape information. Lyttelton et al. \cite{Lyttelton-curvature} computed surface parameterizations that match surface curvature. Fischl et al. \cite{Fischl2} improved the alignment of cortical folding patterns by minimizing the mean squared difference between the average convexity across a set of subjects and that of the individual. Wang et al. \cite{YalinMutual} computed surface registrations that maximize the mutual information between mean curvature and conformal factor maps across subjects. Lord et al. \cite{Lord-isometry} matched surfaces by minimizing the deviation from isometry.

In most situations, one has to pay extra attention to ensure the optimal map computed is diffeomorphic. Hence, developing an effective optimization algorithm that guarantees to give diffeomorphic surface registrations is necessary. This motivates us to look for a simple representation of surface diffeomorphisms which helps to simplify the optimization procedure.

\section{Theoretical Background}
In this section, we describe some basic mathematical concepts related to our algorithms. For details, we refer readers to \cite{Gardiner} and \cite{Lehto}.

A surface $S$ with a conformal structure is called a \emph{Riemann surface}. Given two Riemann surfaces $M$ and $N$, a map $f:M\to N$ is \emph{conformal} if it preserves the surface metric up to a multiplicative factor called the {\it conformal factor}. An immediate consequence is that every conformal map preserves angles. With the angle-preserving property, a conformal map effectively preserves the local geometry of the surface. 
\begin{figure*}[ht]
\centering
\includegraphics[height=2.32in]{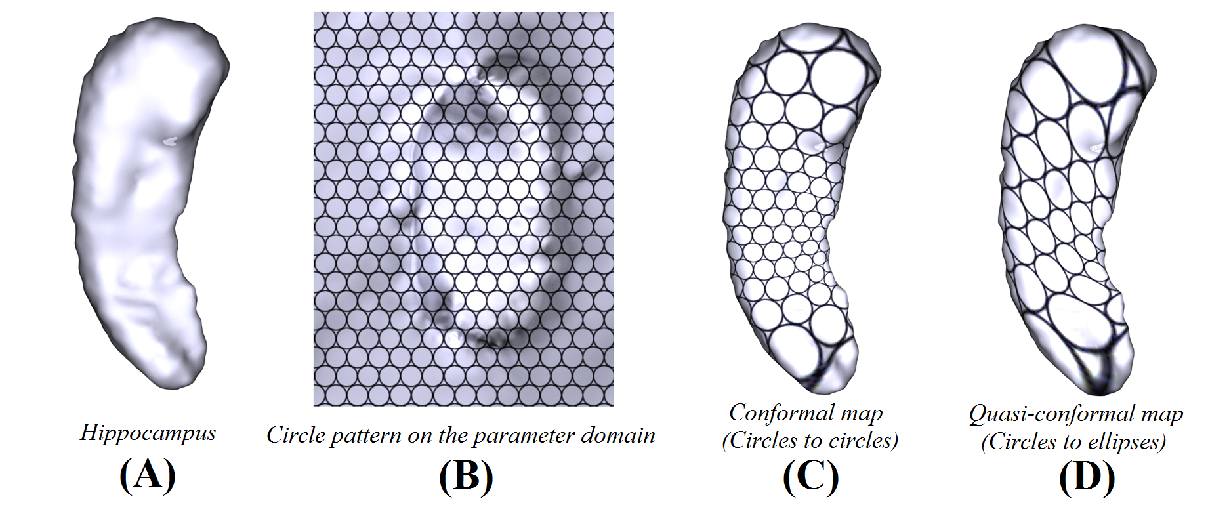}
\caption{Illustration of a conformal map and a quasiconformal map. (A) shows a hippocampal surface. A circle packing pattern is plotted on the parameter domain as in (B). (C) shows a conformal parameterization, which maps circles on the parameter domain to circles on the surface. (D) shows a quasiconformal parameterization, which maps circles on the parameter domain to ellipses on the surface.\label{fig:quasiconformalillustration}}
\end{figure*}
A generalization of conformal maps is \emph{quasi-conformal} maps, which are orientation-preserving diffeomorphisms between Riemann surfaces with bounded conformality distortion, in the sense that their first order approximations takes small circles to small ellipses of bounded eccentricity \cite{Gardiner}. Thus, a conformal homeomorphism that maps a small circle to a small circle can also be regarded as quasi-conformal. Figure \ref{fig:quasiconformalillustration} illustrates the idea of conformal and quasiconformal maps.

Mathematically, $f \colon \mathbb{C} \to \mathbb{C}$ is quasi-conformal provided that it satisfies the Beltrami equation:
\begin{equation}
\frac{\partial f}{\partial \overline{z}} = \mu(z) \frac{\partial f}{\partial z}.
\end{equation}
\noindent for some complex valued functions $\mu$ satisfying $||\mu||_{\infty}< 1$. In terms of the metric tensor, consider the effect of the pullback under $f$ of the usual Euclidean metric $ds_E^2$; the resulting metric is given by:
\begin{equation}
f^*(ds_E^2) = |\frac{\partial f}{\partial z}|^2 |dz + \mu(z) d\overline{z}|^2.
\end{equation}
\noindent which, relative to the background Euclidean metric $dz$ and $d\overline{z}$, has eigenvalues $(1+|\mu|)^2 \left| \frac{\partial f}{\partial z} \right|^2$ and $(1-|\mu|)^2 \left| \frac{\partial f}{\partial z} \right|^2$. $\mu$ is called the \emph{Beltrami coefficient}, which is a measure of non-conformality. In particular, the map $f$ is conformal around a small neighborhood of $p$ when $\mu(p) = 0$. Infinitesimally, around a point $p$, $f$ may be expressed with respect to its local parameter as follows:
\begin{equation}
\begin{split}
f(z) & = f(p) + f_{z}(p)z + f_{\overline{z}}(p)\overline{z} \\
& = f(p) + f_{z}(p)(z + \mu(p)\overline{z}).
\end{split}
\end{equation}
Obviously, $f$ is not conformal if and only if $\mu(p)\neq 0$ at $p$. Inside the local parameter domain, $f$ may be considered as a map composed of a translation to $f(p)$ together with a stretch map $S(z)=z + \mu(p)\overline{z}$, which is postcomposed by a multiplication of $f_z(p)$, which is conformal. All the conformal distortion of $S(z)$ is caused by $\mu(p)$. $S(z)$ is the map that causes $f$ to map a small circle to a small ellipse. From $\mu(p)$, we can determine the angles of the directions of maximal magnification and shrinking and the amount of them as well. Specifically, the angle of maximal magnification is $\arg(\mu(p))/2$ with magnifying factor $1+|\mu(p)|$; The angle of maximal shrinking is the orthogonal angle $(\arg(\mu(p)) -\pi)/2$ with shrinking factor $1-|\mu(p)|$. The distortion or dilation is given by:
\begin{equation}
K = {(1+|\mu(p)|)}/{(1-|\mu(p)|)}.
\end{equation}

\begin{figure*}[ht]
\centering
\includegraphics[height=2.38in]{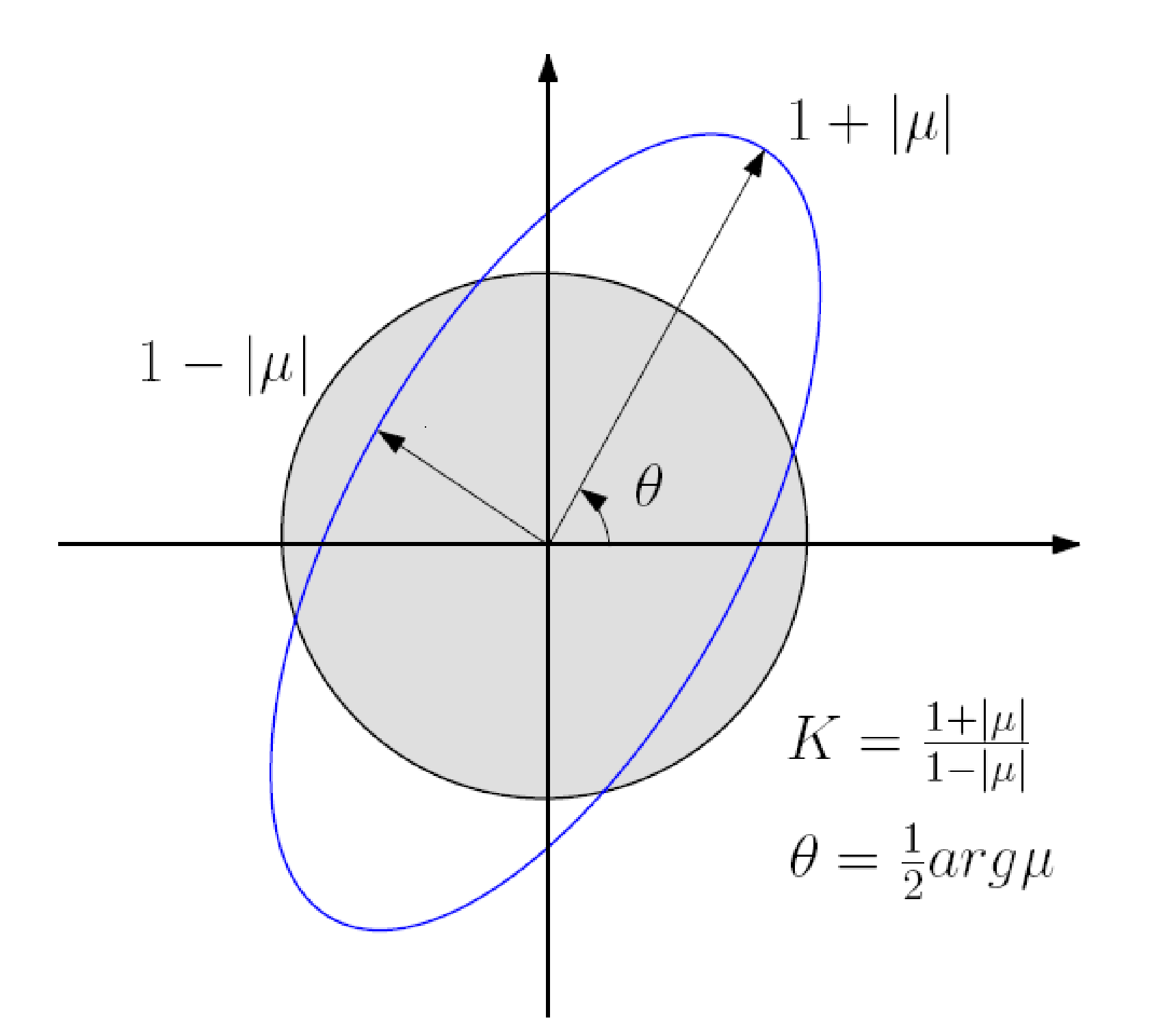}
\caption{Illustration of how the Beltrami coefficient $\mu$ measures the distortion of a quasi-conformal mapping that maps a small circle to an ellipse with dilation $K$.\label{fig:beltramifig}}
\end{figure*}

\noindent Thus, the Beltrami coefficient $\mu$ gives us important information about the properties of the map (See Figure \ref{fig:beltramifig}).

Now, suppose $\mu$ and $\sigma$ are the Beltrami coefficients of the quasiconformal maps $f^{\mu}$ and $f^{\sigma}$ respectively. Then the Beltrami coefficient $\tau$ of the composition map $f^{\tau}=f^{\sigma} \circ (f^{\mu})^{-1}$ can be computed as:
\begin{equation}\label{eqt:composition}
\tau=\left( \frac{\sigma-\mu}{1-\overline{\mu}\sigma} \frac{1}{\theta} \right) \circ (f^\mu)^{-1},
\end{equation}
where $\theta=\frac{\overline{p}}{p}$ and $p=\frac{\partial}{\partial z} f^{\mu}(z)$. In particular, if $f^{\sigma}$ is the identity, that is, if $\sigma=0$, then
\begin{equation}
\tau=-(\mu\frac{p}{\overline{p}}) \circ (f^{\mu})^{-1}.
\end{equation}

\section{Main Algorithm}
In this section, we discuss in detail the main algorithms in this paper. Our goal is to look for a simple representation scheme for the space of surface diffeomorphisms, with the least number of constraints possible, to simplify the optimization process.

\subsection{The Beltrami Holomorphic Flow}
In this part, we describe two theorems about the Beltrami Holomorphic Flow(BHF) on the sphere $\mathbb{S}^2$ and the unit disk $\mathbb{D}$. All the algorithms developed in this paper are mainly based on these theorems.

\begin{thm}[Beltrami holomorphic flow on $\mathbb{S}^2$]\label{thm:BHFC}
There is a one-to-one correspondence between the set of quasiconformal diffeomorphisms of $\mathbb{S}^2$ that fix the points $0$, $1$, and $\infty$ and the set of smooth complex-valued functions $\mu$ on $\mathbb{S}^2$ with $||\mu||_{\infty}=k<1$. Here, we have identified $\mathbb{S}^2$ with the extended complex plane $\overline{\mathbb{C}}$. Furthermore, the solution $f^\mu$ to the Beltrami equation depends holomorphically on $\mu$. Let $\{\mu(t)\}$ be a family of Beltrami coefficients depending on a real or complex parameter $t$. Suppose also that $\mu(t)$ can be written in the form
\begin{equation}
\mu(t)(z)=\mu(z)+t\nu(z)+t\epsilon(t)(z)
\end{equation}
for $z\in\C$, with suitable $\mu$ in the unit ball of $C^\infty(\C)$, $\nu,\epsilon(t)\in L^\infty(\C)$ such that $\parallel\epsilon(t)\parallel_\infty \rightarrow 0$ as $t \rightarrow 0$. Then for all $w \in \C$,
\begin{equation}
f^{\mu(t)}(w)=f^\mu(w)+tV(f^{\mu},\nu)(w)+o(|t|)
\end{equation}
locally uniformly on $\C$ as $t \rightarrow 0$, where
\begin{equation}\label{eqt:flow}
V(f^{\mu},\nu)(w)=-\frac{f^\mu(w)(f^\mu(w)-1)}{\pi}\int_\C \frac{\nu(z)((f^\mu)_z(z))^2}{f^\mu(z)(f^\mu(z)-1)(f^\mu(z)-f^\mu(w))}\,dx\,dy.
\end{equation}
\end{thm}

\bigskip
\begin{proof}
This theorem is due to Bojarski. For detailed proof, please refer to [2].
\end{proof}

Theorem \ref{thm:BHFC} states that any diffeomorphism of $\mathbb{S}^2$ that fixes $0$, $1$ and $\infty$ can be represented uniquely by a Beltrami coefficient. In fact, the 3-point correspondence can be arbitrarily set, instead of fixing $0$, $1$ and $\infty$ only. This can be done easily by composing M\"{o}bius transformations to the diffeomorphism. Let $f:\mathbb{S}^2 \to \mathbb{S}^2$ be any diffeomorphism of $\mathbb{S}^2$. Picking any 3-point coresspondence $\{a, b, c \in \mathbb{S}^2\} \leftrightarrow \{f(a), f(b), f(c) \in \mathbb{S}^2\}$, we can look for unique M\"{o}bius transformations $\phi_1$ and $\phi_2$ that map $\{a,b,c\}$ and $\{f(a),f(b),f(c)\}$ to $0,1,\infty$ respectively. Then, the composition map $\widetilde{f} := \phi_2 \circ f\circ\phi_1^{-1}$ is a diffeomorphism of $\mathbb{S}^2$ that fixes $0$, $1$ and $\infty$ and can be represented by a unique Beltrami coefficient. In other words, given a diffeomorphism $f$ of $\mathbb{S}^2$ and any 3-point correspondence, we can represent $f$ uniquely by a Beltrami coefficient.

The theorem also gives the variation of the diffeomorphism under the variation of the Beltrami coefficient. In order to adjust the diffeomorphism, we can simply adjust the Beltrami coefficient by using the variational formula.

\medskip

Theorem \ref{thm:BHFC} can be further extended to diffeomorphisms of the unit disk $\mathbb{D}$.

\begin{thm}[Beltrami holomorphic flow on $\mathbb{D}$]\label{thm:BHFD}
There is a one-to-one correspondence between the set of quasiconformal diffeomorphisms of $\mathbb{D}$ that fix the points $0$ and $1$ and the set of smooth complex-valued functions $\mu$ on $\mathbb{D}$ for which $||\mu||_{\infty}=k<1$. Furthermore, the solution $f^\mu$ depends holomorphically on $\mu$. Let $\{\mu(t)\}$ be a family of Beltrami coefficients depending on a real or complex parameter $t$. Suppose also that $\mu(t)$ can be written in the form
\begin{equation}
\mu(t)(z)=\mu(z)+t\nu(z)+t\epsilon(t)(z)
\end{equation}
for $z\in\D$, with suitable $\mu$ in the unit ball of $C^\infty(\D)$, $\nu,\epsilon(t)\in L^\infty(\D)$ such that $\parallel\epsilon(t)\parallel_\infty \rightarrow 0$ as $t \rightarrow 0$. Then for all $w \in \D$
\begin{equation}
f^{\mu(t)}(w)=f^\mu(w)+tV(f^{\mu},\nu)(w)+o(|t|)
\end{equation}
locally uniformly on $\D$ as $t \rightarrow 0$, where
\begin{equation}\label{eqt:flow_disk}
\begin{split}
&V(f^{\mu},\nu)(w)=-\frac{f^\mu(w)(f^\mu(w)-1)}{\pi} \\
&\left(\int_\D \frac{\nu(z)((f^\mu)_z(z))^2}{f^\mu(z)(f^\mu(z)-1)(f^\mu(z)-f^\mu(w))}\,dx\,dy +\int_\D \frac{ \overline{\nu(z)}(\overline{(f^\mu)_z(z)})^2} {\overline{f^\mu(z)}(1-\overline{f^\mu(z)})(1-\overline{f^\mu(z)}f^\mu(w))}\,dx\,dy.\right).
\end{split}
\end{equation}
\end{thm}

\bigskip

\begin{proof}
The proof of this theorem can be found in the Appendix.
\end{proof}

Theorem \ref{thm:BHFD} states that any diffeomorphism of $\mathbb{D}$ that fixes 2 points (i.e. $0$ and $1$) can be represented uniquely by a Beltrami coefficient. Again, the 2-point correspondence can be arbitrary. Let $g\colon\mathbb{D}\to \mathbb{D}$ be a diffeomorphism of $\mathbb{D}$. Given any 2-point correspondence $\{a,b\in \mathbb{D}\} \leftrightarrow \{g(a),g(b)\in \mathbb{D}\}$, we can find two unique M\"{o}bius tranformations $\phi_1$ and $\phi_2$ of $\mathbb{D}$ that map $\{a,b\}$ and $\{g(a),g(b)\}$ to $\{0,1\}$ respectively. Then, the composition map $\widetilde{g} := \phi_2 \circ g\circ\phi_1^{-1}$ is a diffeomorphism of $\mathbb{D}$ that fixes $0$ and $1$  and can be represented by a unique Beltrami coefficient. Theorem \ref{thm:BHFD} also gives the variation of the diffeomorphism of $\mathbb{D}$ under the variation of the Beltrami coefficient. Therefore, we can again adjust the diffeomorphism of $\mathbb{D}$ by adjusting the Beltrami coefficient, which is a much simpler functional space.

Theorem \ref{thm:BHFC} and Theorem \ref{thm:BHFD} can be extended to genus 0 closed surfaces and open surfaces with disk topology. Therefore, they can be applied to represent general surface diffeomorphisms. This will be discussed in Section \ref{representation}.

\subsection{Representation of Surface Diffeomorphisms using BCs} \label{representation}
As mentioned earlier, it is crucial to look for a simple representation for the space of all surface diffeomorphisms so that the optimization procedure can be simplified. Surface registration is commonly represented by 3D coordinate functions in $\mathbb{R}^3$. This representation requires lots of storage space and is difficult to manipulate. For example, 3D coordinate functions have to satisfy a constraint in the Jacobian $J$ (namely, $J>0$) in order to preserve the 1-1 correspondence of surface maps. The Jacobian constraint is a complicated partial differential inequality. Enforcing this constraint adds extra difficulty in manipulating and adjusting surface maps. It is therefore important to have a simpler representation with as few constraints as possible.

Theorem \ref{thm:BHFC} and \ref{thm:BHFD} allow us to represent surface diffeomorphisms of $\mathbb{S}^2$ and $\mathbb{D}$ by Beltrami coefficients. The theorems can be further extended to genus 0 closed surfaces and open surfaces with disk topology.

Let $S_1$ and $S_2$ be two genus 0 closed surfaces with a 3-point correspondence between them: $\{p_1, p_2, p_3 \in S_1\} \leftrightarrow \{q_1, q_2, q_3 \in S_2\}$. By Riemann mapping theorem, $S_1$ and $S_2$ can both be uniquely parameterized by conformal maps $\phi_1\colon S_1 \to \mathbb{S}^2$ and $\phi_2\colon S_2 \to \mathbb{S}^2$ respectively, such that $\phi_1(p_1)=0,\phi_1(p_2)=1, \phi_1(p_3)=\infty$ and $\phi_2(q_1)=0,\phi_2(q_2)=1, \phi_2(q_3)=\infty$. The conformal parameterizations can be computed using the discrete Ricci flow method \cite{Jin08TVCGRicci}. Given any surface diffeomorphism $f\colon S_1\to S_2$. The composition map $\widetilde{f}:=\phi_2\circ f\circ\phi_1^{-1}\colon\mathbb{S}^2 \to \mathbb{S}^2$ is a diffeomorphism from $\mathbb{S}^2$ to itself fixing $0$, $1$ and $\infty$. By Theorem \ref{thm:BHFC}, $\widetilde{f}$ can be uniquely represented by a Beltrami coefficient $\widetilde{\mu}$ defined on $\mathbb{S}^2$. Hence, $f$ can be uniquely represented by a Beltrami coefficient $\mu := \widetilde{\mu}\circ \phi_1 ^{-1}$ defined on $S_1$. In other words, we have proven the following:

\bigskip

\begin{cor}\label{Cor1}
Let $S_1$ and $S_2$ be two genus 0 closed surfaces. Suppose $f\colon S_1\to S_2$ is a surface diffeomorphism. Given 3-point correspondence $\{p_1, p_2, p_3 \in S_1\} \leftrightarrow \{f_(p_1), f_(p_2), f_(p_3) \in S_2\}$, $f$ can be represented by a unique Beltrami coefficient $\mu : S_1 \to \mathbb{C}$.
\end{cor}

\bigskip

Similarly, Theorem \ref{thm:BHFD} can be extended to open surfaces with disk topology. Let $M_1$ and $M_2$ be two genus 0 open surfaces. Given two points correspondence $\{p_1, p_2\in M_1\} \leftrightarrow \{q_1, q_2\in M_2\}$ between them. We can again uniquely parameterize $M_1$ and $M_2$ conformally to map the corresponding points to $0$ and $1$. Denote them by $\phi_1\colon M_1 \to \mathbb{D}$ and $\phi_2\colon M_2 \to \mathbb{D}$. The composition map $\widetilde{f}:=\phi_2\circ f\circ\phi_1^{-1}\colon \mathbb{D} \to \mathbb{D}$ is a diffeomorphism of $\mathbb{D}$ fixing $0$ and $1$. Again, $\widetilde{f}$ can be uniquely represented by a Beltrami coefficient $\widetilde{\mu}$ defined on $\mathbb{D}$. Hence, $f$ can be uniquely represented by a Beltrami coefficient $\mu := \widetilde{\mu}\circ \phi_1 ^{-1}$ defined on $M_1$. So, we have the following Corollary:

\bigskip

\begin{cor}\label{Cor2}
Let $M_1$ and $M_2$ be two genus 0 open surfaces with disk topology. Suppose $f\colon M_1\to M_2$ is a surface diffeomorphism. Given 2-point correspondence $\{p_1, p_2\in M_1\} \leftrightarrow \{f_(p_1), f_(p_2)\in M_2\}$, $f$ can be represented by a unique Beltrami coefficient $\mu\colon M_1 \to \mathbb{C}$.
\end{cor}

\bigskip

Corollary \ref{Cor1} and \ref{Cor2} allows us to represent diffeomorphisms of genus 0 closed surfaces and open surfaces with disk topology using Beltrami coefficients. Thus, we can use the Beltrami coefficient $\mu_f$ associated uniquely to such diffeomorphism $f$ to represent $f$. First of all, we need to compute the Beltrami coefficient $\widetilde{\mu}_{\widetilde{f}}$ of the composition map $\widetilde{f} = \phi_2\circ f\circ\phi_1^{-1}\colon D\to D$, where $D$ is the common conformal parameter domain of the surfaces.. Mathematically, $\widetilde{\mu}_{\widetilde{f}}$ is given by the following formula:
\begin{equation}\label{eqt:mutilde}
\begin{split}
\widetilde{\mu}_{\widetilde{f}} &= \frac{\partial \widetilde{f}}{\partial \overline{z}}/ \frac{\partial \widetilde{f}}{\partial z}= \frac{1}{2}(\frac{\partial \widetilde{f}}{\partial x} + \sqrt{-1}\frac{\partial \widetilde{f}}{\partial y})/\frac{1}{2}(\frac{\partial \widetilde{f}}{\partial x} - \sqrt{-1}\frac{\partial \widetilde{f}}{\partial y}).
\end{split}
\end{equation}

Then, the Beltrami coefficient $\mu_f$ can be computed by $\mu_f := \widetilde{\mu}_{\widetilde{f}}\circ \phi_1^{-1}\colon S_1 \to \mathbb{C}$. $\mu_f$ is a complex-valued functions defined on $S_1$ with $||\mu_f||_{\infty} <1$. There are no restrictions on $\mu_f$ that it has to be 1-1, surjective or satisfy some constraints on the Jacobian. With this representation, we can easily manipulate and adjust surface maps without worrying about destroying their diffeomorphic property.

In practice, surfaces are commonly approximated by discrete meshes comprising of triangular or rectangular faces. They are parameterized onto the mesh $D$ in $\mathbb{C}$. Then the partial derivatives (or gradient) of the map can be discretely approximated on each face of $D$. By taking average, the partial derivatives and hence the Beltrami coefficient can be computed on each vertex. The detailed numerical implementation can be found in the Appendix.

Besides adjusting surface maps preserving diffeomorphism, another advantage of Beltrami coefficients is that they consist of two real functions only, namely the real and imaginary parts. Compared to the representation of surface maps using 3D coordinate functions, this representation reduces 1/3 of the original storage space.

The computational algorithm can be summarized as follows:

\noindent $\mathbf{Algorithm\ 1.}$  {\it Beltrami Representation of Surface Diffeomorphisms}\\
\noindent {\it Input: Surface diffeomorphism $f\colon S_1\to S_2$; points correspondence $\{p_i\} \leftrightarrow \{q_i = f(p_i)\}$.}\\
\noindent {\it Output: Beltrami representation $\mu_f\colon S_1 \to \mathbb{C}$ of $f\colon S_1\to S_2$.}
\begin{enumerate}
\item {\it Compute the conformal parameterizations of $S_1$ and $S_2$ that map $\{p_i\}$ and $\{q_i\}$ to consistent locations on the parameter domain $D$. Denote them by $\phi_1\colon S_1 \to D$ and $\phi_2\colon S_2 \to D$}
\item {\it Set $\widetilde{f} = \phi_2\circ f\circ\phi_1^{-1}\colon D\to D$ and compute the Beltrami coefficient $\widetilde{\mu}_{\widetilde{f}}$ by Equation \ref{eqt:mutilde}.}
\item {\it Compute the Beltrami coefficient $\mu_f\colon S_1 \to \mathbb{C}$ using $\mu_f := \widetilde{\mu}_{\widetilde{f}}\circ \phi_1^{-1}$}.
\end{enumerate}
\begin{figure*}[ht]
\centering
\includegraphics[height=3in]{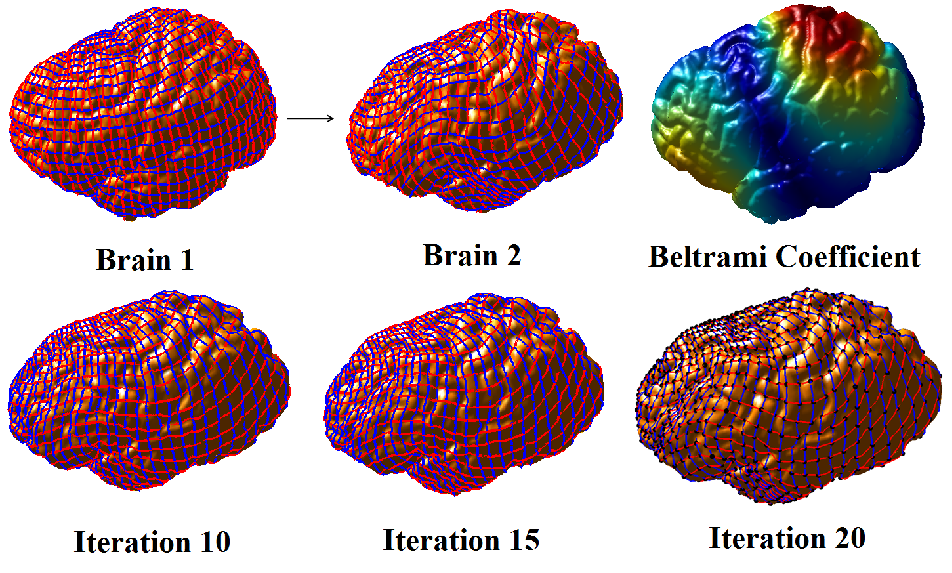}
\caption{Beltrami representation and reconstruction of a surface diffeomorphism $f$ on the brain surface. The top-left figure shows a surface diffeomorphism between two different brain surfaces. The top-right figure shows the Beltrami representation $\mu$ of $f$. Here, the colormap of $|\mu|$ is shown. The bottom row shows the reconstructed map after different number of iterations $N$ using BHF reconstruction. When $N=20$, the map closely resembles the original map (The black dots show the exact positions under the original map.) \label{fig:brainreconstruction}}
\end{figure*}

\subsection{Reconstruction of Surface Diffeomorphisms from BCs}
Given the Beltrami coefficient $\mu$ defined on $S_1$. It is important to have a reconstruction scheme to compute the associated quasi-conformal diffeomorphism $f^{\mu}$. This allows us to move back and forth between BCs and surface diffeomorphisms. We propose the {\it Beltrami holomorphic flow} (BHF) method to reconstruct the surface diffeomorphism $f^{\mu}\colon S_1\to S_2$ associated with a given $\mu$. BHF iteratively flows the identity map to $f^{\mu}$. In this part, we describe the BHF reconstruction method in detail.

\begin{figure*}[ht]
\centering
\includegraphics[height=3.65in]{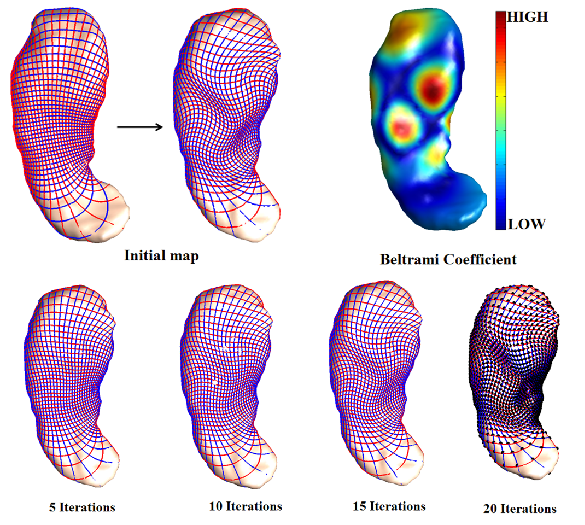}
\caption{Beltrami representation and reconstruction of a surface diffeomorphism $f$ on hippocampal surfaces. The top-left figure shows a surface diffeomorphism between two different hippocampal surfaces. The top-right figure shows the Beltrami representation $\mu$ of $f$. The colormap of $|\mu|$ is shown. The bottom row shows the reconstructed map after different number of iterations $N$ using BHF reconstruction. When $N=20$, the map closely resembles the original map (The black dots shows the exact positions under the original map.)\label{fig:hipporeconstruction}}
\end{figure*}
The variation of $f^{\mu}$ under the variation of $\mu$ can be expressed explicitly. Suppose $\widetilde{\mu}(z) = \mu(z) + t\nu (z) +o(|t|)$ where $z = x + iy \in \mathbb{C}$. Then, $f^{\widetilde{\mu}(z)}(w) = f^{\mu}(w) + tV({f}^{\mu},\nu)(w)  +o(|t|)$, where
\begin{equation} \label{BHFiteration}
\begin{split}
V({f}^{\mu},\nu)&(w) =\int_D K(z,w)\, dx\,dy,
\end{split}
\end{equation}

\noindent and
\begin{equation}\label{BHFiteration2}
K(z,w) =
 \begin{cases}
 -\frac{f^\mu(w)(f^\mu(w)-1)}{\pi}\left(\frac{\nu(z)((f^\mu)_z(z))^2}{f^\mu(z)(f^\mu(z)-1)(f^\mu(z)-f^\mu(w))}\right) & D = \mathbb{S}^2,\\
 -\frac{f^\mu(w)(f^\mu(w)-1)}{\pi}\left(\frac{\nu(z)((f^\mu)_z(z))^2}{f^\mu(z)(f^\mu(z)-1)(f^\mu(z)-f^\mu(w))}+ \frac{\overline{\nu(z)}(\overline{(f^\mu)_z(z)})^2}{\overline{f^\mu(z)}(1-\overline{f^\mu(z)})(1-\overline{f^\mu(z)}f^\mu(w))}\right) & D = \mathbb{D}.
 \end{cases}
\end{equation}

We can also write $V({f}^{\mu},\nu)(w)$ as:
\begin{equation}\label{BHFiteration3}
\begin{split}
V({f}^{\mu},\nu)&(w) =\int_D  \left( \begin{array}{c}
G_1 \nu_1 + G_2 \nu_2\\
G_3 \nu_1 + G_4 \nu_2\end{array} \right) \,dx\,dy,
\end{split}
\end{equation}

\noindent where $\nu = \nu_1 + i \nu_2$ and $G_1$, $G_2$, $G_3$, $G_4$ are real-valued functions defined on $D$. Here, we identify $A+iB$ as $\left( \begin{array}{c}
A\\
B\end{array} \right)$.

Using this fact, we propose to use BHF to iteratively flow the identity map to $f^{\mu}$. Given the parameterizations $\phi_1\colon S_1 \to D$ and $\phi_2\colon S_2 \to D$, we look for the map $\widetilde{f}^{\mu} = \phi_2 \circ f^{\mu}\circ\phi_1^{-1}\colon D\to D$ associated uniquely with $\widetilde{\mu}= \mu\circ\phi_1^{-1}\colon D\to \mathbb{C}$. $f^{\mu}$ can then be obtained by $f^{\mu} = \phi_2^{-1} \circ \widetilde{f}^{\mu}\circ\phi_1$.

We start with the identity map $\mathbf{Id}$ of which the Beltrami coefficient is identically equal to 0. Let $N$ be the number of iterations. Define $\widetilde{\mu}_k = k\widetilde{\mu}/N$, $k=\{0,1,2,\ldots N\}$. Let $\widetilde{f}^{\widetilde{\mu}_k }$ be the map associated with $\widetilde{\mu}_k$. Note that $\widetilde{f}^{\widetilde{\mu}_0} = \mathbf{Id}$ and $\widetilde{f}^{\widetilde{\mu}_N} = \widetilde{f}^{\widetilde{\mu}}$. Equation \ref{BHFiteration} allows us to iteratively compute $\widetilde{f}^{\widetilde{\mu}_k}$ and thus obtain a sequence of maps flowing from $\mathbf{Id}$ to $\widetilde{f}^{\widetilde{\mu}}$. The iterative scheme is given by:
\begin{equation}\label{iterationscheme}
\begin{split}
& \widetilde{f}^{\widetilde{\mu}_{k+1}} = \widetilde{f}^{\widetilde{\mu}_{k}} + V(\widetilde{f}^{\widetilde{\mu}_{k}},\frac{\widetilde{\mu}}{N});\ \ \ \widetilde{f}^{\widetilde{\mu}_{0}} =\mathbf{Id}
\end{split}
\end{equation}

The computational algorithm of the reconstruction scheme can be summarized in Algorithm 2. The detailed numerical implementation can be found in Appendix.

\medskip

\noindent $\mathbf{Algorithm\ 2.}$  {\it Reconstruction of Surface Diffeomorphisms from BCs}\\
\noindent {\it Input: Beltrami Coefficient $\mu$ on $S_1$; conformal parameterizations of $S_1$ and $S_2$: $\phi_1$ and $\phi_2$; Number of iterations $N$}\\
\noindent {\it Output: Surface diffeomorphism $f^{\mu}\colon S_1\to S_2$ associated to $\mu$.}
\begin{enumerate}
\item {\it Set $k=0$; $\widetilde{f}^{\widetilde{\mu}_0} = \mathbf{Id}$.}
\item {\it Set $\widetilde{\mu}_k := k\widetilde{\mu}/N$; Compute $\widetilde{f}^{\widetilde{\mu}_{k+1}} = \widetilde{f}^{\widetilde{\mu}_{k}} + V(\widetilde{f}^{\widetilde{\mu}_{k}},\frac{\widetilde{\mu}}{N})$; $k = k+1$.}
\item {\it Repeat Step 2 until $k=N$; Set $f^{\mu} := \phi_2^{-1} \circ \widetilde{f}^{\widetilde{\mu}}\circ\phi_1\colon S_1 \to S_2$}.
\end{enumerate}

\medskip

\begin{figure*}[ht]
\centering
\includegraphics[height=2.65in]{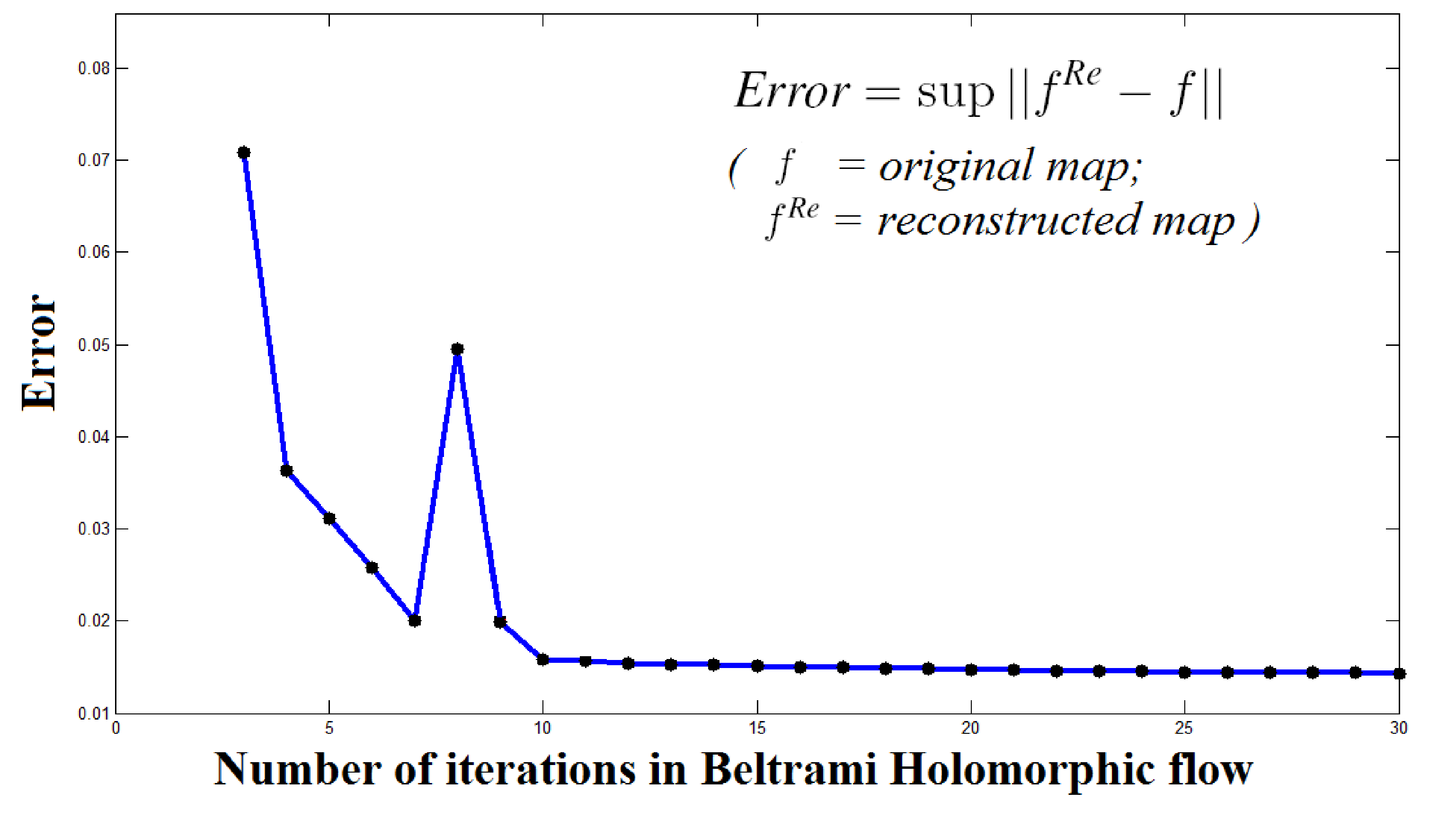}
\caption{The error of the reconstructed map $f^{\mathrm{Re}}$ versus the number of iterations used in the BHF process. \label{fig:reconstructederror}}
\end{figure*}

Figure \ref{fig:brainreconstruction} and \ref{fig:hipporeconstruction} illustrate the idea of reconstructing surface diffeomorphisms from BCs on human brain surfaces and hippocampal surfaces respectively. BHF computes a sequence of surface maps $\{\widetilde{f}^{\widetilde{\mu}_k}\}$ converging to $\widetilde{f}^{\widetilde{\mu}}$. The approximation of $\widetilde{f}^{\widetilde{\mu}_k}$ is more accurate with a smaller time step, or equivalently, a larger number of iterations $N$. Figure \ref{fig:reconstructederror} shows the error of the reconstructed map $f^{\mathrm{Re}}$ versus different number of iterations $N$ used in the BHF process. The error is defined as $Error = \sup ||f^{\mathrm{Re}} - f||$, where $f$ is the original map. As expected, the error decreases as $N$ increases. In practice, the approximations are very accurate when $N\geq 15$. In our experiments, we set $N=20$.

\subsection{BHF Optimization of Surface Registrations}
We have described a simple representation scheme for surface diffeomorphisms using BCs.  The space of BCs is a simple functional space with the least amount of constraints. There are no restrictions requiring BCs to be 1-1, surjective or satisfy some constraints on its Jacobian. With BCs, we can easily manipulate and adjust surface maps, while ensuring the diffeomorphic(1-1, onto and smooth) property of the surface registration.

Theorem \ref{thm:BHFC} and \ref{thm:BHFD} give us the variation of surface maps under the variation of their BCs (Equation \ref{BHFiteration} and \ref{BHFiteration2}). This allows us to perform optimization on the space of BCs, instead of working directly on the space of surface diffeomorphisms. The diffeomorphic property of the optimal surface registration can also be easily ensured during the optimization process.

Given an energy functional $E$ defined on the space of surface diffeomorphisms, we can easily reformulate $E$ and redefine it on the space of BCs. With the BHF variation, we can derive the Euler-Lagrange equation on $E$ to optimize BCs iteratively. To demonstrate the idea, we consider a simple example of optimizing surface maps between two human brain surfaces.
\begin{figure*}[ht]
\centering
\includegraphics[height=1.75in]{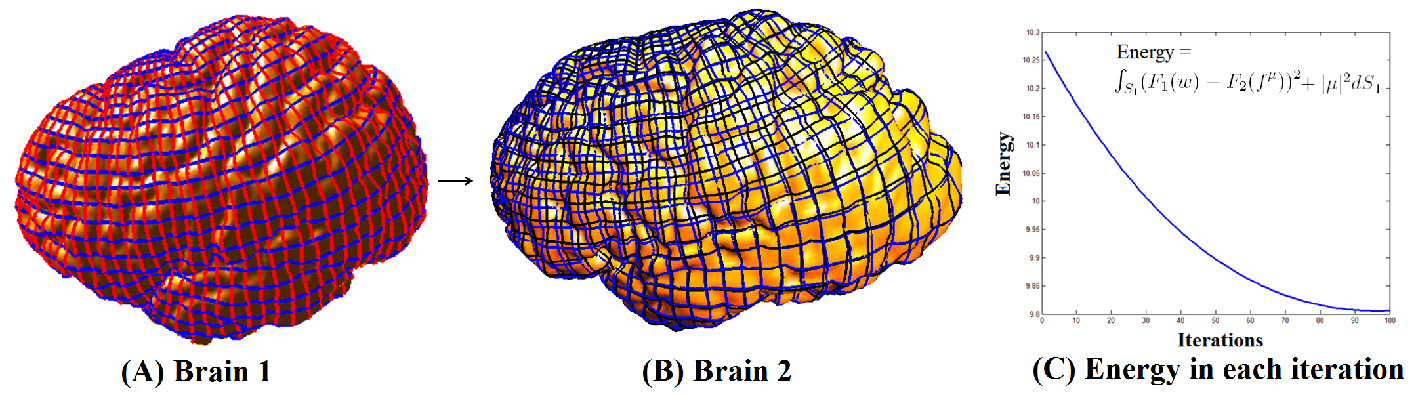}
\caption{Illustration of BHF optimization scheme on brain surfaces. This example shows the optimization result of matching two feature functions $F_1$ and $F_2$ on the two brain surfaces. The blue grid represents the initial map, while the black grid represents the optimized map. \label{fig:example2_2}}
\end{figure*}
\paragraph*{Example 4.1} Consider two different human brain surfaces $S_1$ and $S_2$ as shown in Figure \ref{fig:example2_2}. Denote the conformal parameterizations of them by $\phi_1\colon S_1 \to D$ and $\phi_2\colon S_2 \to D$. In surface registration, it is often important to find an optimal 1-1 correspondence that matches some intensity feature functions defined on each surfaces. Let $F_1\colon S_1 \to \mathbb{R}$ and $F_2: S_2 \to \mathbb{R}$ be two intensities (functions) defined on $S_1$ and $S_2$ respectively. As an illustration, we define $F_1$ and $F_2$ as $F_1 := \phi_1^{-1}(5.2x^2 + 3.3y^2)$ and $F_2 := \phi_2^{-1}(6.8x^2 + 2.8y)$. We propose to find $f:S_1 \to S_2$ minimizing $E(f) = \int_{S_1} (F_1(w) - F_2(f(w)))^2 +|\mu(w)|^2 \,dw$. The optimized map $f$ is a quasi-conformal map that best matches $F_1$ and $F_2$ while preserving the conformality as good as possible. We can formulate the energy functional to be defined on the space of BCs over the conformal parameter domain $D$. That is,
\begin{equation}
E(\mu) = \int_D (F_1(w) - F_2(f^{\mu}))^2 +|\mu(w)|^2 \,dw
\end{equation}

The Euler-Lagrange equation can be derived as follow:
\begin{equation}
\begin{split}
\frac{d}{dt}|_{t=0} E(\mu + t\nu) &= \int_D \frac{d}{dt}|_{t=0} \left( (F_1(w) - F_2(f^{\mu+t\nu}(w)))^2 +|\mu(w) + t\nu(w)|^2 \right)\,dw \\
&= -\int_D 2(F_1 - F_2(f^{\mu}))\nabla F_2(f^{\mu}) \frac{d}{dt}|_{t=0} f^{\mu + t\nu} - 2\mu\cdot \nu \,dw\\
&= -\int_D \int_D \left( \begin{array}{c}
A\\
B\end{array} \right)\cdot \left( \begin{array}{c}
G_1\nu_1 +G_2 \nu_2\\
G_3\nu_1 +G_4 \nu_2\end{array} \right) \,dz - 2\mu\cdot \nu \,dw\\
&= -\int_D \left( \int_D \left( \begin{array}{c}
AG_1 +BG_3\\
AG_2 +BG_4 \end{array} \right) \,dw - \left( \begin{array}{c}
2\mu_1\\
2\mu_2 \end{array} \right)\right)\cdot \left( \begin{array}{c}
\nu_1\\
\nu_2 \end{array}\right) \,dz,
\end{split}
\end{equation}
\noindent where $\left( \begin{array}{c} A\\B\end{array} \right)= 2(F_1- F_2(f^{\mu}))\nabla F_2$; $\mu = \mu_1 + i\mu_2$ and $\nu = \nu_1 + i\nu_2$.

So, the descent direction for $\mu = \mu_1 + i\mu_2$ is
\begin{equation}\label{eqt:descent}
\frac{d\mu_1}{dt} = \int_D (AG_1 + BG_3)\,dw - 2\mu_1 \mathrm{\ and\ } \frac{d\mu_2}{dt} = \int_D (AG_2 + BG_4)\,dw - 2\mu_2.;
\end{equation}

We can iteratively optimize the energy $E$ as follow:
\begin{equation}
\mu^{n+1}=\mu^n+dt \left( \begin{array}{c}
\int_D (A^n G_1^n + B^n G_3^n)\,dw - 2\mu_1\\
\int_D (A^n G_2^n + B^n G_4^n)\,dw - 2\mu_2\end{array} \right)
\end{equation}

Figure \ref{fig:example2_2} shows the experimental result for this example. (A) shows the standard grid on Brain 1. The standard grid is mapped by the initial map to Brain 2, which is shown as the blue grid. We optimize the map such that it minimizes the energy functional. The resulting map is plotted as the black grid. (C) shows the energy at each iteration. It decreases as the number of iterations increases. This shows that our BHF optimization algorithm can iteratively optimize the energy functional.
\hspace*{\fill}~\QED\par\endtrivlist\unskip

\bigskip
Therefore, with BHF, we can perform optimizations over the space of BCs, which is a much simpler functional space with least amount of contraints, and simplify the optimization procedure significantly.

\section{Applications}
In this section, we outline some applications of our proposed optimization algorithm to surface registration. These applications are motivated from practical problems encountered in medical imaging.

\subsection{Optimized Conformal Parameterization with Landmark Matching}
With BHF, we first develop an algorithm to effectively compute {\it landmark-matching optimized conformal maps} between surfaces. A landmark-matching optimized conformal map refers to a map that matches corresponding landmarks across surfaces, while preserving conformality as much as possible. It is very important for research applications in computational anatomy. For example, in human brain mapping, neuroscientists are often interested in finding a 1-1 correspondence between brain surfaces that matches sulcal/gyral landmark curves, which are important anatomical features \cite{AutoLui}. Besides matching these brain features, they also want the maps to preserve local geometry as much as possible. Conformal maps are best known to preserve local geometry and hence are commonly used. However, landmark matching cannot be guaranteed under conformal maps. Therefore, it is of interest to look for maps which are as conformal as possible and match landmarks well.

Most existing algorithms for computing landmark-matching optimized conformal maps cannot ensure exact landmark matching. Some existing algorithms can align landmarks consistently, but bijectivity is usually not guaranteed especially when a large number of landmark constraints are imposed \cite{Lui07}. Here, we introduce a variational approach to compute an optimized conformal map iteratively by minimizing the $L^2$ norm of a Beltrami coefficient $\mu$ (the Beltrami energy). Since $\mu$ is a measure of local distortion in conformality, our proposed algorithm is in fact looking for the best landmark-matching map, which is as conformal as possible. Also, a map is bijective as long as $||\mu||_{\infty} <1$. This can be easily controlled and guaranteed in each iteration by minimizing the Beltrami energy in our algorithm.

\begin{figure*}[ht]
\centering
\includegraphics[height=2.5in]{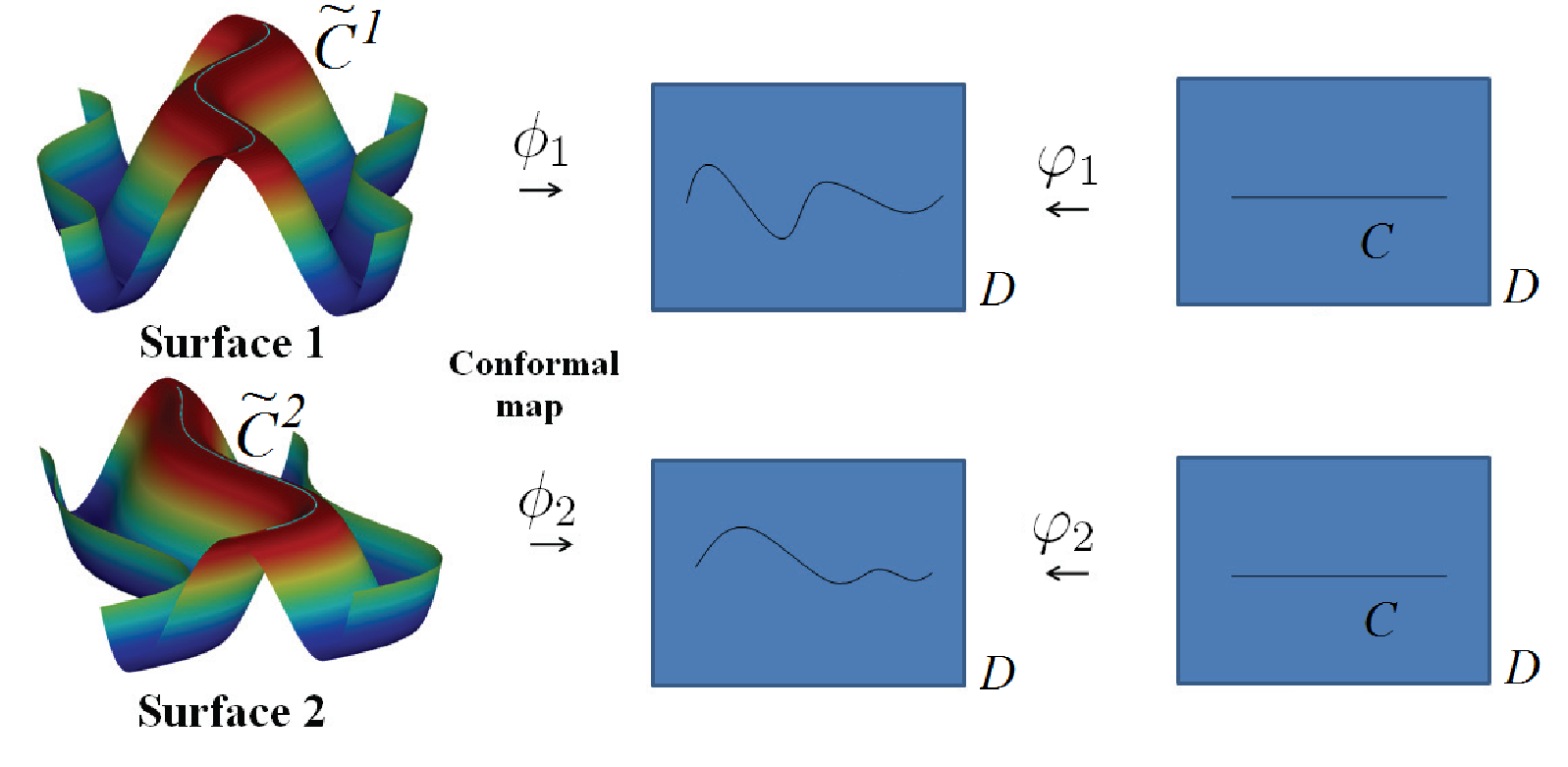}
\caption{This figure shows the framework of the landmark-matching optimized conformal parameterization algorithm. \label{fig:landmarkillustration}}
\end{figure*}

Given two surfaces $S_1$ and $S_2$ with the same topology. Denote the corresponding landmark curves on $S_1$ and $S_2$ by $\{\widetilde{C}_k^1\}$ and $\{\widetilde{C}_k^2\}$ respectively. We first parameterize $S_1$ and $S_2$ conformally onto a common parameter domain $D$ ($=\mathbb{D}$ or $\mathbb{S}^2 \cong \overline{\mathbb{C}}$). Let $\phi_1\colon S_1 \to D$ and $\phi_2\colon S_2 \to D$ be the parameterizations. We proposed to look for two maps $\varphi_1\colon D\to D$ and $\varphi_2\colon D\to D$ such that $\varphi_i^{-1}$ ( $i=1,2$) maps landmarks $\{\phi_i(\widetilde{C}_k^i)\}$ onto the consistent landmarks $\{C_k\}$ on $D$ (see Figure \ref{fig:landmarkillustration}), and that it minimizes the following energy functional:
\begin{equation} \label{eqt:landmarkenergy}
E(\varphi_i) = \int_D |\mu_{\varphi_i}|^2.
\end{equation}

\begin{figure*}[ht]
\centering
\includegraphics[height=2.3in]{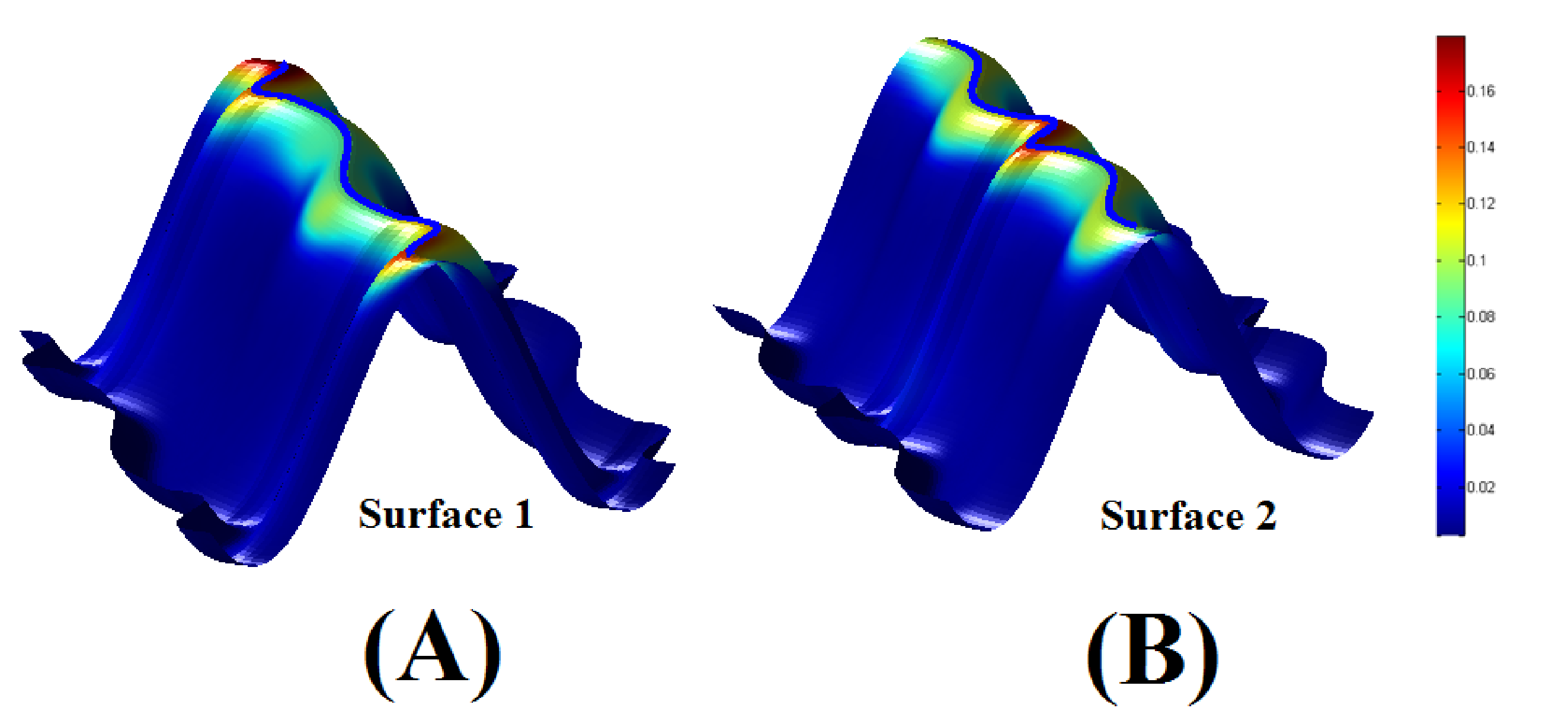}
\caption{The Beltrami coefficients of the optimized conformal parameterizations of 2 synthetic surfaces fixing 1 landmark. The norms of the Beltrami coefficients are plotted as colormap, which are very small except near the landmark curves. It means that the conformality distortion is accumulated around the landmarks.\label{fig:beltrami12}}
\end{figure*}

Equation \ref{eqt:landmarkenergy} ensures that each landmark-matching parametrization $\varphi_i$ has the least conformality distortion. Hence, the local geometric distortion under $\varphi_i$ is minimized. Starting from the conformal map with $\mu=0$, the energy also ensures the property that $|\mu|_{\infty} < 1$ and so the diffeomorphic property of the minimizer is guaranteed. A landmark-matching map $f$ between $S_1$ and $S_2$ can then be obtained by the composition map: $f :=\phi_1^{-1}\circ\varphi_2\circ\varphi_1^{-1}\circ\phi_i$. We can compute the Euler-Lagrange equation of Equation \ref{eqt:landmarkenergy} with respect to $\mu_{\varphi_i}$ as follow:
\begin{equation} \label{E-LEnergylandmark}
\begin{split}
\frac{d}{dt}|_{t=0} E(\mu_{\varphi_i} + tv) &= \int_D \frac{d}{dt}|{t=0} |\mu_{\varphi_i} + tv|^2 \\
&= 2\int_D [\mathbf{Re}(\mu_{\varphi_i})\mathbf{Re}(v)) + \mathbf{Im}(\mu_{\varphi_i})\mathbf{Im}(v))]
\end{split}
\end{equation}
\begin{figure*}[Ht!]
\centering
\includegraphics[height=3.96in]{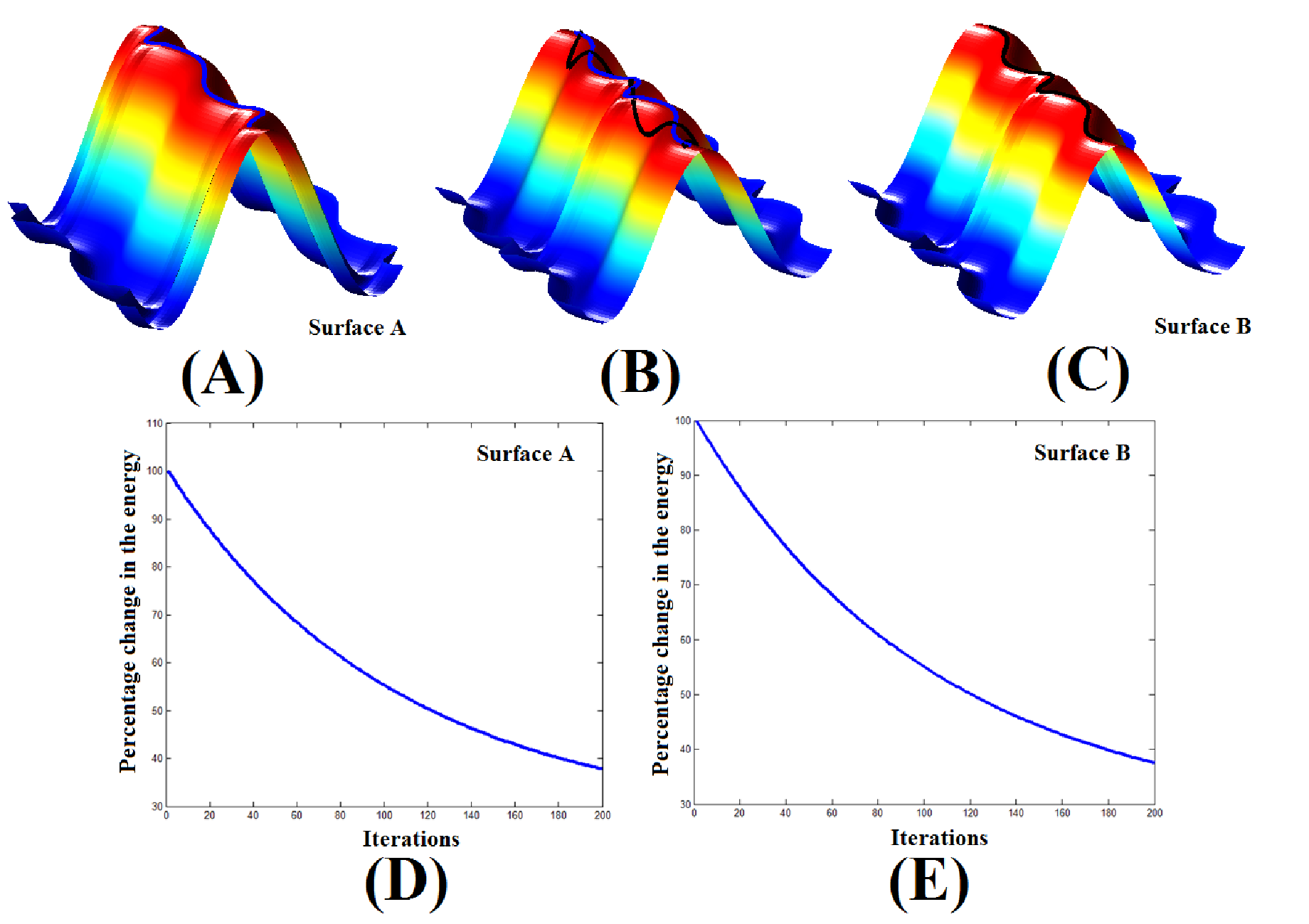}
\caption{Illustration of landmark-matching optimized conformal parameterizations of synthetic surfaces with 1 landmark. The blue curves in (A) and (B) represent the landmarks on the two surfaces. The landmark on surface A cannot be mapped onto the one on surface B under the conformal map(the black curve in (B)). With optimized conformal parameterization, the corresponding landmarks on each surface can be exactly matched, as shown in (C). The percentage change in energy functionals of the optimized conformal parameterizations for surface A and B are shown in (D) and (E).\label{fig:landmarksyn}}
\end{figure*}
The derivative in Equation \ref{E-LEnergylandmark} is negative when $v = - 2\mu_{\varphi_i}$. Hence, we can iteratively minimize $E(\mu_{\varphi_i})$ by the following scheme:
\begin{equation}
\mu_{\varphi_i}^{n+1} - \mu_{\varphi_i}^{n} = - 2\mu_{\varphi_i}^n dt.
\end{equation}
The detailed computational algorithm can be described as follow:

\noindent $\mathbf{Algorithm\ 3.}$  {\it Optimized Conformal Parameterization with Landmark Matching}\\
\noindent {\it Input: Surfaces $S_1$ and $S_2$; Landmark curves $\widetilde{C}_k^1$ on $S_1$, $\widetilde{C}_k^2$ on $S_2$.}\\
\noindent {\it Output: Optimized conformal parameterization $\varphi_1$ and $\varphi_2$ of $S_1$ and $S_2$ with landmark matching.}
\begin{enumerate}
\item {\it Compute the initial map $\varphi_i^0$ that aligns landmark curves $\{\phi_i(\widetilde{C}_k^i)\}$ to $\{C_k\}$ on $D$. Set $n=0$.}
\item {\it Compute the Beltrami coefficient $\mu_{\varphi_i}^{n}$ of $\varphi_i^n$. Let $\mu_{\varphi_i}^{n+1} = \mu_{\varphi_i}^{n} - 2\mu_{\varphi_i}^n dt$.}
\item {\it Compute $\vec{V}_n= V(\varphi_i^n, -2\mu_{\varphi_i}^n)$ using the BHF formula}.
\item {\it Let $\varphi_i^{n+1}(p) = \varphi_i^n(p) + \delta(p)\vec{V}_n(p)dt$ , where $\delta$ is a smooth delta function on $D$ that is equal to zero around $\{C_k^i\}$ and one elsewhere. This ensures landmarks are matched in each iteration. Set n = n+1.}
\item {\it Repeat Step 2 to Step 5. If $|E(\mu_{\varphi_i}^{n+1})-E(\mu_{\varphi_i}^{n})| <\epsilon$, $\mathbf{stop}$.}
\end{enumerate}
We have tested our proposed method on synthetic data as well as real medical data. Figure \ref{fig:landmarksyn} shows the result of matching two synthetic surfaces with one landmark on each surface. The blue curves on (A) and (B) represent the landmarks on the two surfaces. Under a conformal map, the landmark on surface A cannot be mapped exactly onto the one on surface B (the black curve in (B)). Using our proposed method, the corresponding landmarks on each surface can be exactly matched, as shown in (C). (D) and (E) show the percentage change in energy functionals of the optimized conformal parameterizations for surface A and B. The energies decrease as the number of iterations increases. This indicates a decrease in the conformalilty distortion. Figure \ref{fig:beltrami12} shows the Beltrami coefficient of each optimized conformal parameterization. The colormap shows the norm of the Beltrami coefficient. Note that the norm of the Beltrami coefficient is very small except near the landmark curve. It means the conformality distortion is accumulated around the landmarks as expected.

\begin{figure*}[ht]
\centering
\includegraphics[height=4.2in]{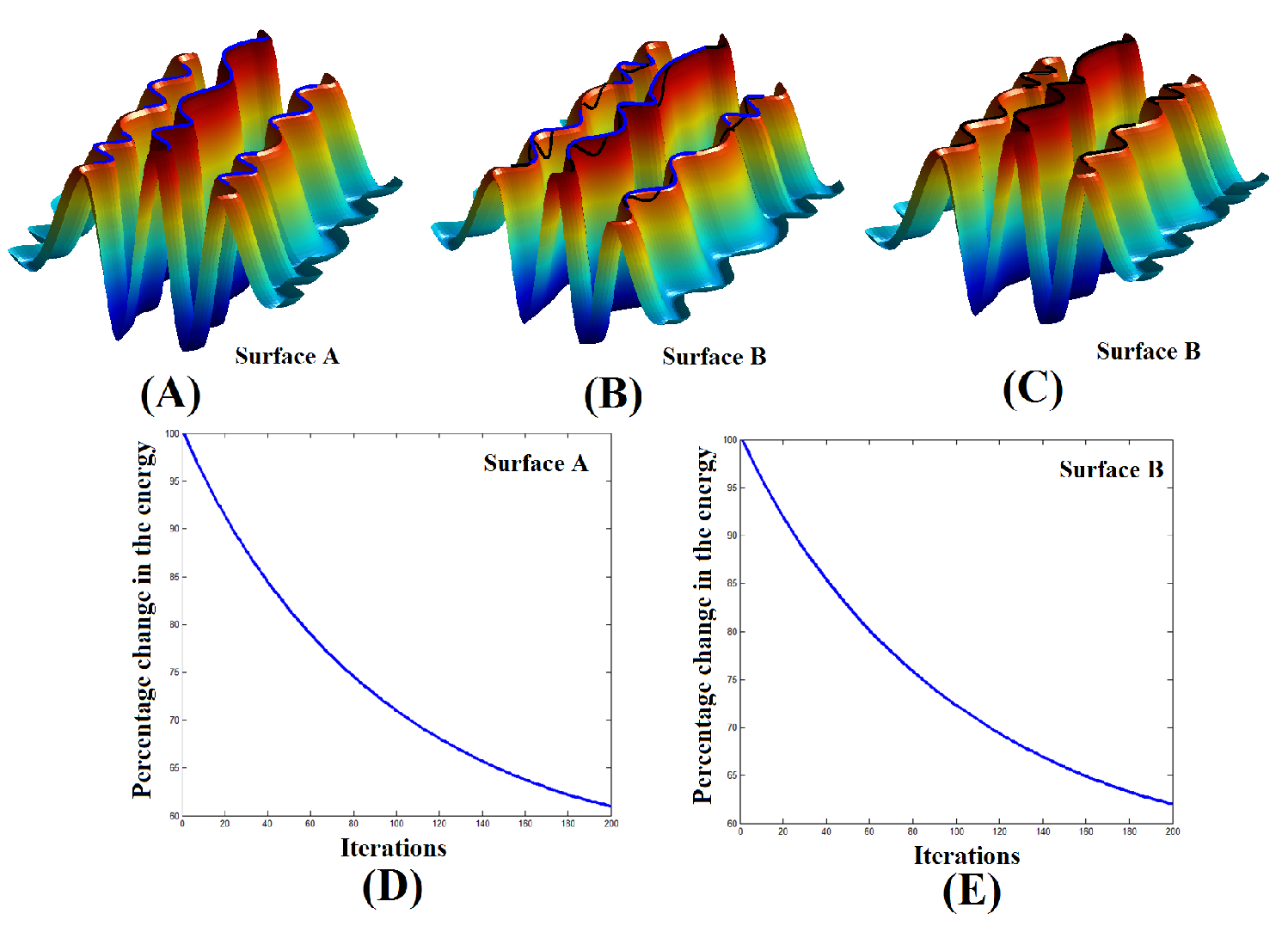}
\caption{Landmark-matching optimized conformal parameterizations of 2 synthetic surfaces fixing 5 landmarks. The blue curves on (A) and (B) represent the landmarks on the surfaces. The landmark on Surface A cannot be mapped to the landmark on Surface B  under the conformal map(black curves in (B)). With optimized conformal parameterization, the corresponding landmarks on each surface can be exactly matched (shown in (C)). The percentage change in energies of the optimized conformal parameterizations for surface A and surface B are plotted in (D) and (E).\label{fig:example3.1}}
\end{figure*}

\begin{figure*}[ht]
\centering
\includegraphics[height=2.1in]{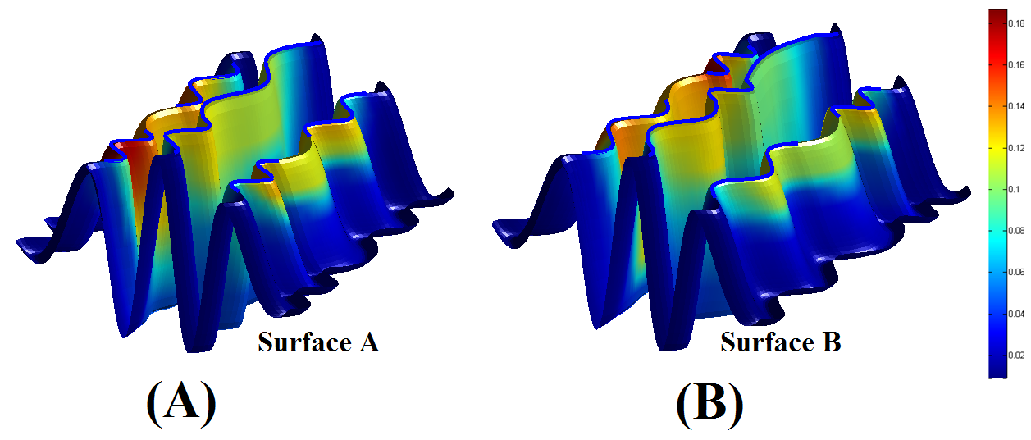}
\caption{The Beltrami coefficient of the optimized conformal parameterization of 2 synthetic surfaces fixing 5 landmarks. The norm of the Beltrami coefficient is very small except near the landmark curve. It means that the conformality distortion is accumulated around the landmarks.\label{fig:example3.2}}
\end{figure*}

We have also tested our algorithm on synthetic surfaces with five landmarks as shown in Figure \ref{fig:example3.1}. Again, the landmarks cannot be exactly matched under a conformal map (see black curves in (B)). However, they are exactly matched using our proposed algorithm. As shown in (D) and (E), the percentage change in energies decreases as the number of iterations increases, meaning that conformality distortion is progressively reduced. Figure \ref{fig:example3.2} shows the Beltrami coefficients of the optimized conformal parameterizations fixing landmarks. Again, the norm of the Beltrami coefficient is very small except near the landmark curves.


\begin{figure*}[ht]
\centering
\includegraphics[height=4.2in]{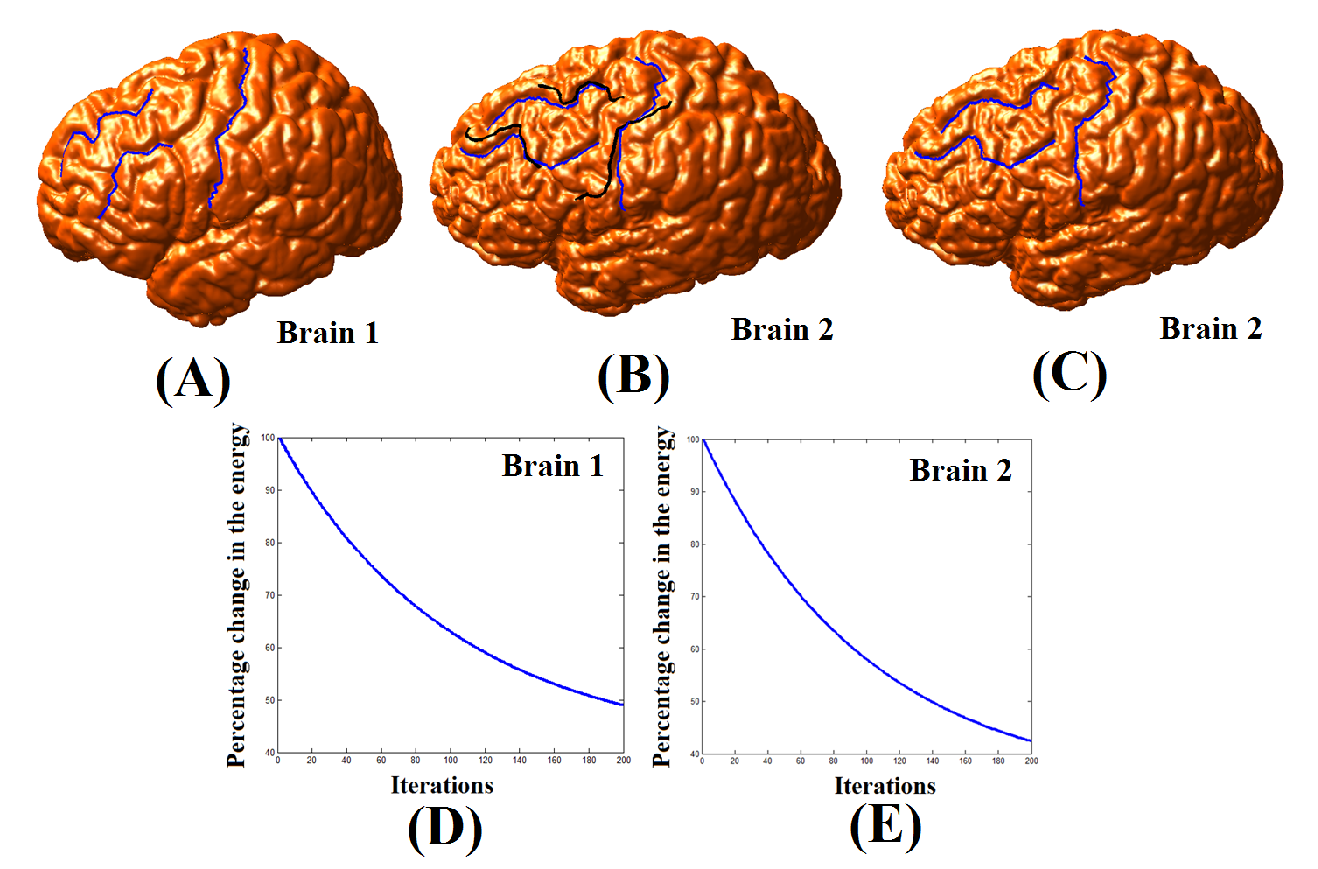}
\caption{Landmark-matching optimized conformal parameterizations of cortical hemispheric surfaces with 3 major sulcal landmarks. The blue curves on (A) and (B) represent the landmarks on the two surfaces. Under a conformal map, the landmarks on Brain A cannot be correctly mapped onto landmarks on Brain B (black curves in (B)). With landmark-matching optimized conformal parameterization, the corresponding landmarks on each surface can be exactly matched as shown in (C). The percentage change in energies of the optimized conformal parameterizations for Brain A and Brain B are shown in (D) and (E).\label{fig:brainlandmark2}}
\end{figure*}

\begin{figure*}[ht]
\centering
\includegraphics[height=2.3in]{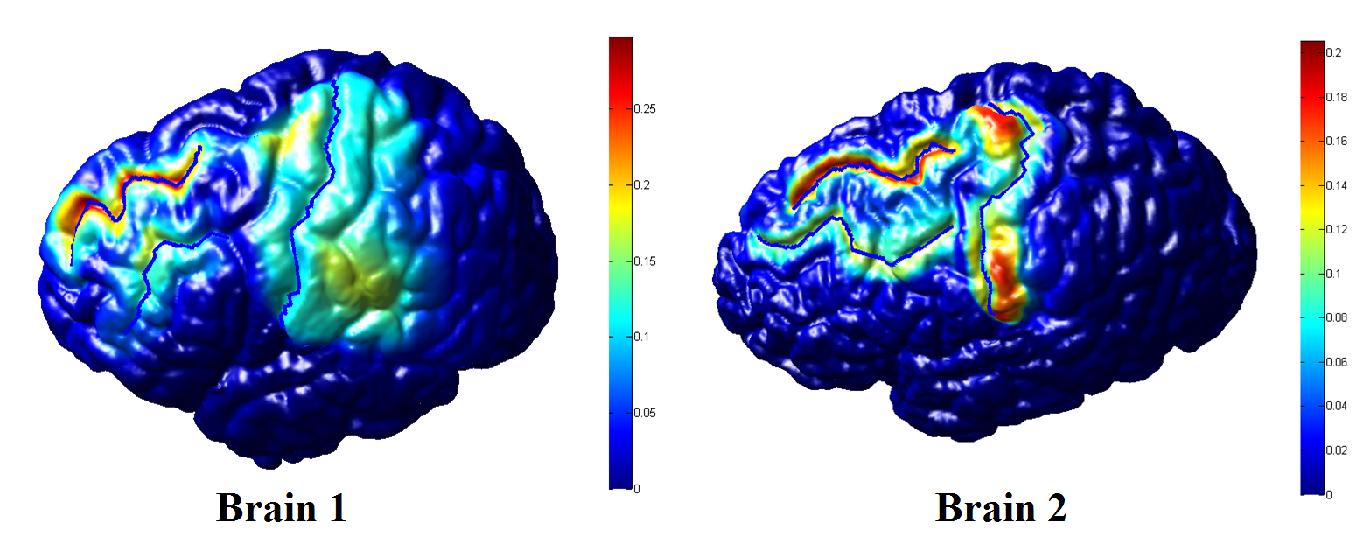}
\caption{The Beltrami coefficient of the optimized conformal parameterizations of 2 cortical hemispheric surfaces with 3 major sulcal landmarks. The norms of the Beltrami coefficients are shown as colormap, which are very small except near the landmark curves. As expected, the conformality distortion is accumulated around the landmarks.\label{fig:brainbel}}
\end{figure*}

Finally, we have tested our algorithm on real cortical hemispheric surfaces extracted from brain MRI scans, acquired from normal subjects at 1.5 T(on a GE Signa scanner). Figures \ref{fig:brainlandmark2}(A) and (B) show 2 different brain surfaces with 3 major sulcal curves labeled on each of them (see the blue curves). Under a conformal map, landmarks on Brain 1 and Brain 2 are not exactly matched (see the black curves in (B)). They are, however, exactly matched using our proposed algorithm as shown in (C). (D) and (E) show the percentage change in the energies of the optimized conformal parameterizations of the surfaces. The energies decrease as the number of iterations increases. This shows that the conformality distortion is gradually reduced. Figure \ref{fig:brainbel} shows the Beltrami coefficients of the optimized conformal parameterizations of the 2 brain surfaces. Again, the norms of the Beltrami coefficients are very small except near the sulci curves.

\begin{figure*}[ht]
\centering
\includegraphics[height=3.85in]{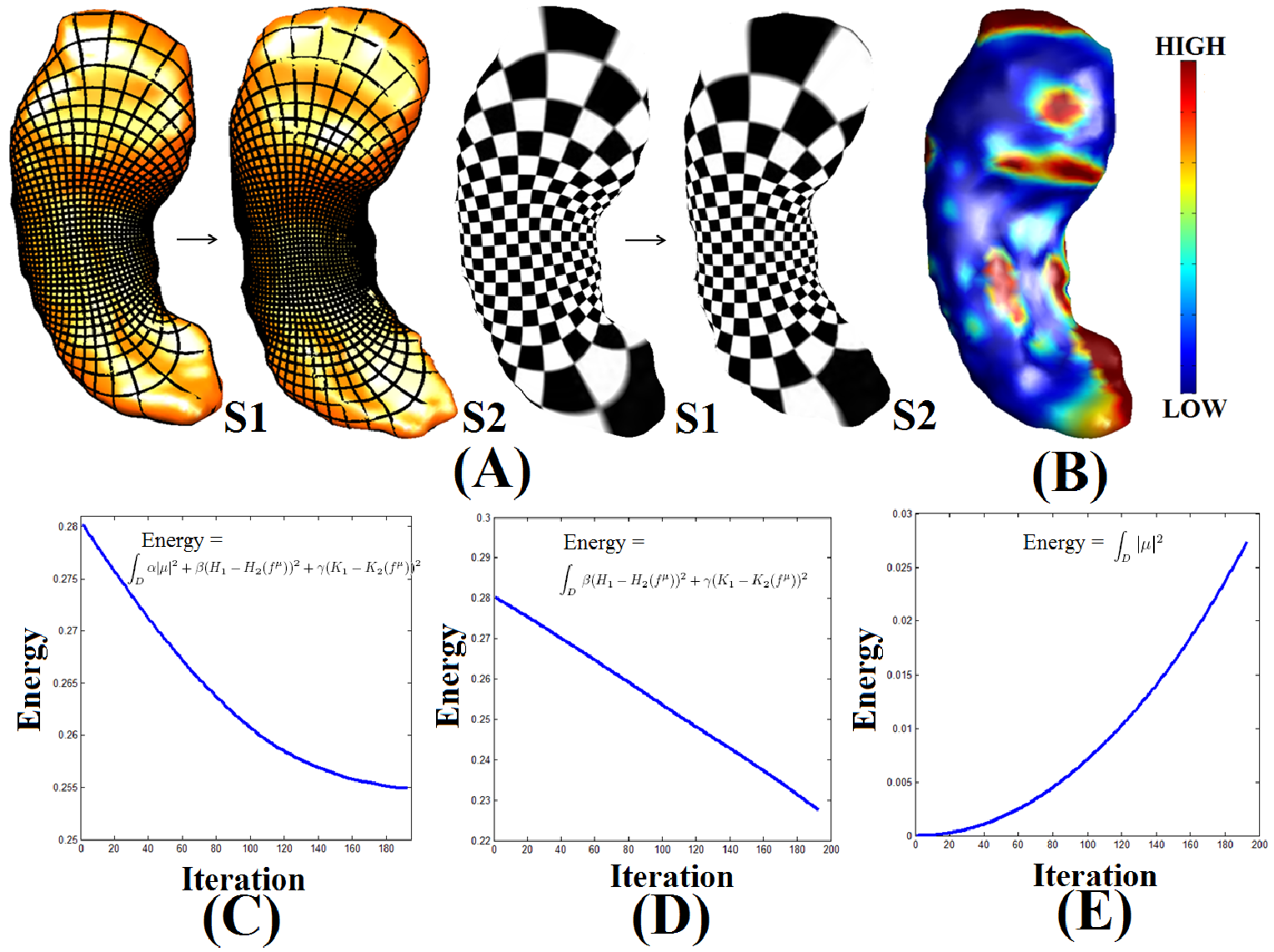}
\caption{Shape registration with geometric matching using Beltrami Holomorphic Flow (BHF). The registration is visualized as grid map and texture map as shown in (A). The optimal shape energy is shown in (B). The percentage changes of the shape energy, curvature mismatch energy and Beltrami energy after different number of iteration are shown in (D), (E) and (F) respectively.\label{fig:hippo3}}
\end{figure*}

\subsection{Hippocampal Registration with Geometric Matching}
In medical imaging, there are cases where anatomical landmark features cannot be easily defined on some brain structures. In such cases, landmark-matching constraints cannot be used as a criterion to establish good correspondence between surfaces. Finding the best 1-1 correspondence between these structures becomes challenging. One typical example is the hippocampus(HP), which is an important structure in the human brains. It belongs to the limbic system and plays important roles in long-term memory and spatial navigation. Surface-based shape analysis is commonly used to study local changes of HP surfaces due to pathologies such as Alzheimer disease (AD), schizophrenia and epilepsy \cite{Thompsonnew}. On HP surfaces, there are no well-defined anatomical landmark features. High-field structural or functional imaging, where discrete cellular fields are evident \cite{Zeinehnew}, is still not routinely used. Finding meaningful registrations between HP surfaces becomes challenging. It is thus important to develop methods to look for good registrations between different HP surfaces without landmarks. To achieve this, we develop an algorithm to automatically register HP surfaces with complete geometric matching and avoid the need to manually label landmark features. This is done by optimizing a compounded energy, which minimizes the $L^2$ norm of the Beltrami coefficient and matches curvatures defined on each surface. Given two hippocampal surfaces $S_1$ and $S_2$. The compounded energy $E_{shape}$ is defined mathematically as
\begin{equation}\label{eqt:shapeenergy}
E_{\mathrm{shape}}(\mu) = \alpha\int_{D} |\mu|^2 + \beta\int_{D} (H_1 - H_2(f^{\mu}))^2 + \gamma\int_{D} (K_1 - K_2(f^{\mu}))^2
\end{equation}
\noindent where $H_1$, $H_2$ are the mean curvatures on $S_1$, $S_2$ respectively, defined on the common parameter domain $D$, and $K_1,\ K_2$ are the Gaussian curvatures. The first integral minimizes the conformality distortion of the surface registration, and ensures the diffeomorphic property of the minimizer by controlling $\mu$. The second and third integrals ensure the optimized registration matches the curvatures as much as possible. It turns out that $E_{\mathrm{shape}}$ is a complete shape index which measures the dissimilarity between two surfaces. Specifically, $E_{\mathrm{shape}} = 0$ if and only if $S_1$ and $S_2$ are geometrically equal up to a rigid motion. Therefore, surface map minimizing $E_{\mathrm{shape}}$ is the best registration that matches the geometric information as much as possible.
\begin{figure*}[ht]
\centering
\includegraphics[height=2in]{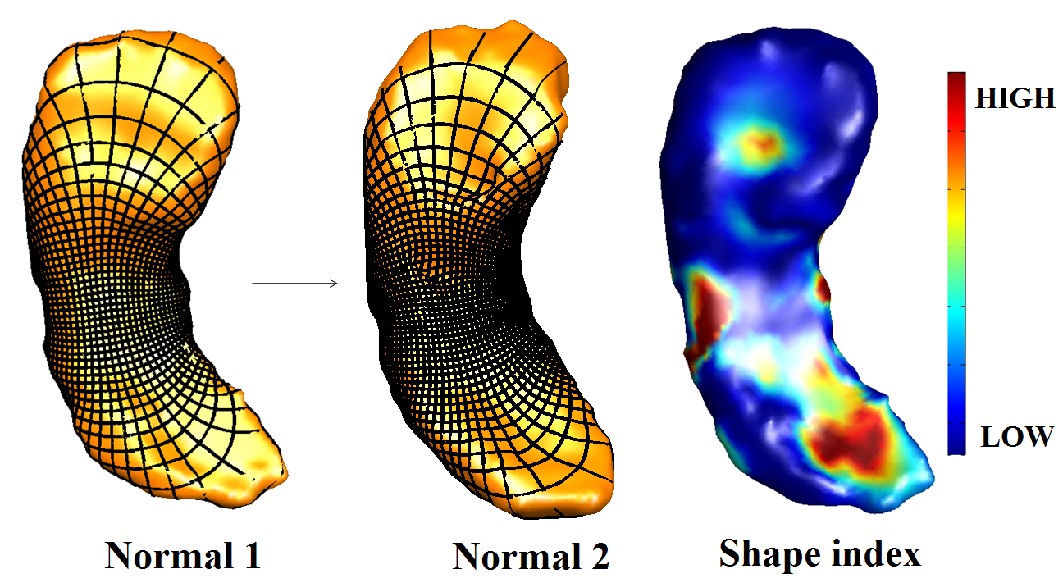}
\caption{BHF registration between two normal subjects. The shape index $E_{\mathrm{shape}}$ is plotted on the right, which captures local shape differences.\label{fig:ex8}}
\end{figure*}
\begin{figure*}[ht]
\centering
\includegraphics[height=2in]{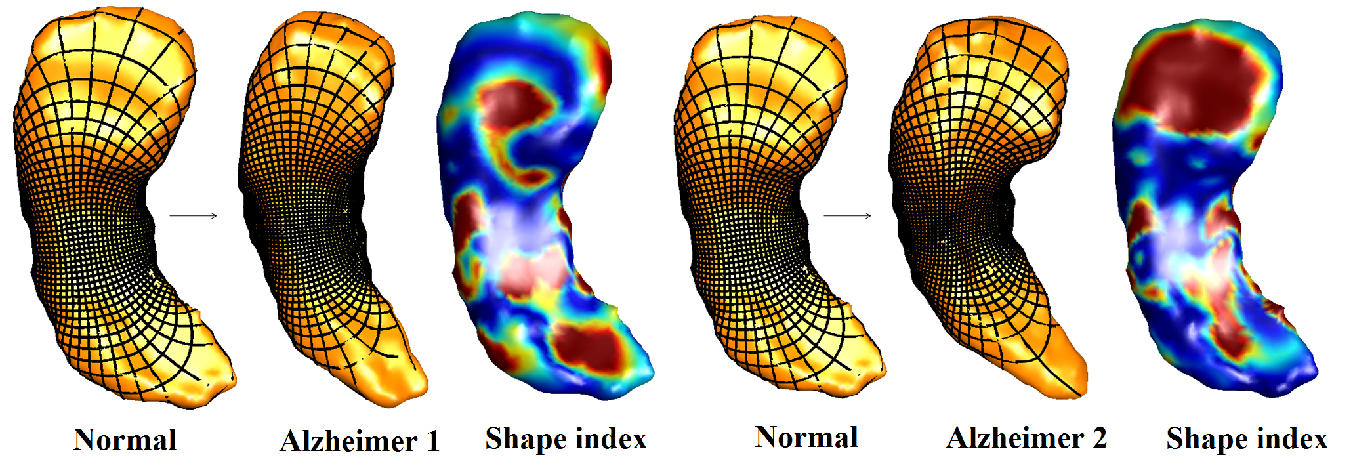}
\caption{BHF registration between 2 normal subjects and 2 subjects with Alzheimer's disease. The local shape differences captured by $E_{\mathrm{shape}}$ are plotted on the surfaces.\label{fig:combined2}}
\end{figure*}
\begin{figure*}[ht]
\centering
\includegraphics[height=2.5in]{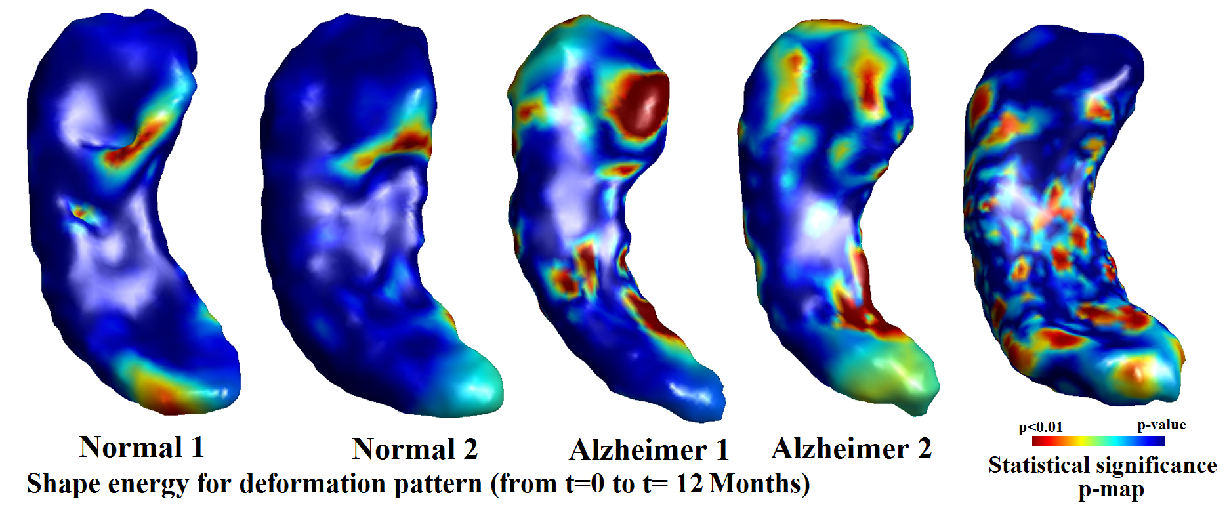}
\caption{Temporal HC shape changes of normal and subjects with Alzheimer's disease.\label{fig:combinedtimepmap2}}
\end{figure*}
We can minimize $E_{\mathrm{shape}}$ in Equation \ref{eqt:shapeenergy} iteratively, using the proposed BHF optimization algorithm. The Euler-Lagrange equation of Equation \ref{eqt:shapeenergy} can be computed as follows:
\begin{equation}\label{E-LEnergygeometric}
\begin{split}
\frac{d}{dt}|_{t=0}E_{\mathrm{shape}} (\mu) &= \alpha\int_D \frac{d}{dt}|_{t=0}|\mu + tv|^2 + \beta\int_D \frac{d}{dt}|_{t=0}(H_1 - H_2(f^{\mu+tv}))^2\\
& + \gamma\int_D \frac{d}{dt}|_{t=0}(K_1 - K_2(f^{\mu+tv}))^2\\
& = 2\alpha\int_D \mu \cdot v - 2\beta\int_D (H_1 - H_2(f^{\mu}))\nabla H_2(f^{\mu}) \cdot \frac{df^{\mu+tv}}{dt}|_{t=0} \\
&\ \  - 2\gamma\int_D (K_1 - K_2(f^{\mu}))\nabla K_2(f^{\mu}) \cdot \frac{df^{\mu+tv}}{dt}|_{t=0}\\
& = 2\int_w \{\alpha\mu(w) - \int_z [(\beta\widetilde{H} +\gamma\widetilde{K})\cdot  \left( \begin{array}{c}
G_1\\
G_2 \end{array} \right), (\beta\widetilde{H} + \gamma\widetilde{K})\cdot \left( \begin{array}{c}
G_3\\
G_4 \end{array} \right) ]\} \cdot v(w)
\end{split}
\end{equation}

\noindent where $\int_w \bullet := \int_D \bullet\ \, dw$ and $\int_z \bullet := \int_D \bullet\ \, dz$ are defined as the integral over variables $w$ and $z$ respectively, $\widetilde{H} := (H_1 - H_2(f^{\mu}))\nabla H_2(f^{\mu})$, $\widetilde{K} := (K_1 - K_2(f^{\mu}))\nabla K_2(f^{\mu})$, $G_i$ is as defined in Equation \ref{BHFiteration3}.

The derivative in Equation \ref{E-LEnergygeometric} is negative when $v = - 2(\mu(w) - \int_z [(\widetilde{H} +\widetilde{K})\cdot G, \mathbf{det}(\widetilde{H} + \widetilde{K}, G)]\ )$. Hence, we can iteratively minimize $E(\mu)$ by the following iterative scheme:
\begin{equation}\label{BHFgeometry}
\mu^{n+1} - \mu^{n} = -2(\alpha\mu^n - \int_z [(\beta\widetilde{H}^n +\gamma\widetilde{K}^n)\cdot G^n, \mathbf{det}(\beta\widetilde{H}^n + \gamma\widetilde{K}^n, G^n)]\ ) \, dt
\end{equation}

The detailed computational algorithm can be described as follows:

\noindent $\mathbf{Algorithm\ 4.}$  {\it BHF Registration with Geometric Matching}\\
\noindent {\it Input: Hippocampal surfaces $S_1$ and $S_2$, step length $dt$, threshold $\epsilon$}\\
\noindent {\it Output: Geometric matching registration $f^{\mu}$ and the shape index $E(f^{\mu})$}
\begin{enumerate}
\item {\it Compute the conformal parameterizations of $S_1$ and $S_2$. Denote them by $\phi_1\colon S_1 \to D$ and $\phi_2: S_2 \to D$}
\item {\it Set $\varphi^0 := \mathbf{Id}\colon D\to D$ and $n=0$.}
\item {\it Compute the Beltrami coefficient $\mu_{\varphi}^{n}$ of $\varphi ^n$ (e.g. $\mu_{\varphi}^{0} = 0$). Update $\mu_{\varphi}^{n+1}$ by Equation \ref{BHFgeometry}}.
\item {\it Compute: $\vec{V}_n= V(\varphi ^n, \mu_{\varphi}^{n+1} - \mu_{\varphi}^{n})$ using Equation \ref{BHFiteration}. Let $\varphi^{n+1} = \varphi^n +\vec{V}_n$. Set n = n+1.}
\item {\it Repeat Step 3 to Step 5. If $|E(\mu_{\varphi}^{n+1})-E(\mu_{\varphi}^{n})| <\epsilon$, $\mathbf{Stop}$.}
\end{enumerate}

We have tested our algorithm on 212 HP surfaces automatically extracted from 3D brain MRI scans with a validated algorithm \cite{Morranew}. Scans were acquired from normal and diseased (AD) elderly subjects at 1.5 T (on a GE Signa scanner). In our experiments, we set $\alpha=1$ and $\beta=\gamma=2$. Experimental results show that our proposed algorithm is effective in registering HP surfaces with geometric matching. The proposed shape energy can also be used to measure local shape difference between HPs. Figure \ref{fig:hippo3}(A) shows two different HP surfaces. They are registered using our proposed BHF algorithm with geometric matching. The registration is visualized using a grid map and a texture map, which shows a smooth 1-1 correspondence. The optimal shape index $E_{\mathrm{shape}}$ is plotted as a colormap in (B). $E_{\mathrm{shape}}$ effectively captures the local shape difference between the surfaces. (C) shows the shape energy in each iteration. With the BHF algorithm, the shape energy decreases as the number of iterations increases. (D) shows the curvature mismatch energy ($E = \int \beta(H_1 - H_2(f))^2 + \gamma(K_1 - K_2(f))^2$). It decreases as the number of iterations increases, meaning that the geometric matching improves. (E) shows the Beltrami coefficient of the map in each iteration, which shows the conformality distortion of the map. Some conformality is intentionally lost to allow better geometric matching.

Figure \ref{fig:ex8} shows the BHF registration between two normal HPs. The complete shape index $E_{\mathrm{shape}}$ is plotted as colormap on the right. Again, $E_{\mathrm{shape}}$ can accurately capture local shape differences between the normal HP surfaces.

Figure \ref{fig:combined2} shows the BHF hippocampal registrations between normal elderly subjects and subjects with Alzheimer's disease. The BHF registrations give smooth 1-1 correspondence between the HP surfaces. We can use the complete shape index $E_{\mathrm{shape}}$ to detect local shape differences between healthy and unhealthy subjects.

We also study the temporal shape changes of normal and AD HP surfaces, as shown in Figure \ref{fig:combinedtimepmap2}. We compute the deformation pattern of its HP surfaces for each subject, measured at time = 0 and time = 12 months (see \cite{Morranew2} for longitudinal scanning details). The left two panels show the temporal deformation patterns for two normal subjects. The middle two panels show the temporal deformation patterns for two AD subjects. The last column shows the statistical significance p-map measuring the difference in the deformation pattern between the normal (n=47) and AD (n=53) groups, plotted on a control HP. The deep red color highlights regions of significant statistical difference. This method can potentially be used to study factors that influence brain changes in AD.

\section{Conclusion}
In this paper, we propose a simple representation of bijective surface maps using Beltrami coefficients(BCs), which helps the optimization process of surface registrations. To complete the representation scheme, we develop a reconstruction algorithm of the surface diffeomorphism from a given BC using the Beltrami holomorphic flow method. This allows us to move back and forth between BCs and surface diffeomorphisms. By formulating the variation of the associated surface map under the variation of BC, we reformulate variational problems over the space of surface diffeomorphisms into variational problems over the space of BCs. It greatly simplifies the optimization procedure. More importantly, a bijective surface map is always guaranteed during the optimization process. Experimental results on synthetic examples and real medical applications show the effectiveness of our proposed algorithms for surface registration.

\bigskip

\section{Appendix}

\begin{center}
\underline{$\mathbf{I.\ Numerical\ Implementation}$}
\end{center}
In this part, we give detailed numerical implementation on how the proposed algorithms can be computed. In practice, all surfaces are represented by meshes which consist of vertices, edges, and triangular/rectangular faces. In our iterative scheme, the functions and their partial derivatives are defined on each vertex and then linearly interpolated to define their values inside each triangular/rectangular face.

\noindent {\it 1. Computation of the Beltrami Coefficient}

Let $f=(f_1,f_2)$ be the diffeomorphism defined on the parameter domain $D$. The Beltrami coefficient $\mu_f$ associated uniquely to $f$ can be computed as follows (see Equation \ref{eqt:mutilde}):
\begin{equation}
\mu_f = [(\frac{\partial f_1}{\partial x} - \frac{\partial f_2}{\partial y}) + i(\frac{\partial f_2}{\partial x} + \frac{\partial f_1}{\partial y})]/[(\frac{\partial f_1}{\partial x} + \frac{\partial f_2}{\partial y}) + i(\frac{\partial f_2}{\partial x} - \frac{\partial f_1}{\partial y})].
\end{equation}

In order to compute $\mu_f$, we simply need to approximate the partial derivatives at each vertex: $D_x f_i(\vec{v}) \approx \frac{\partial f_i}{\partial x}(\vec{v})$ and  $D_y f_i(\vec{v}) \approx \frac{\partial f_i}{\partial y}(\vec{v})$. We first approximate the gradient $\nabla _T f_i$ on each face $T$ by solving:
\begin{equation}
\left( \begin{array}{c}
\vec{v}_1 - \vec{v}_0\\
\vec{v}_2 - \vec{v}_0\end{array} \right)\nabla_T f_i = \left( \begin{array}{c}
\frac{f_i(\vec{v}_1) - f_i(\vec{v}_0)}{|\vec{v}_1 - \vec{v}_0|}\\
\frac{f_i(\vec{v}_2) - f_i(\vec{v}_0)}{|\vec{v}_2 - \vec{v}_0|}\end{array} \right),
\end{equation}
\noindent where $[\vec{v_0},\vec{v_1}]$ and $[\vec{v_0},\vec{v_2}]$ are two edges on $T$. After the gradient $\nabla _T f_i$ have been computed for each face $T$, $D_x f_i(\vec{v})$ and $D_y f_i(\vec{v})$ can be computed by taking average as follows:
\begin{equation}
\left( \begin{array}{c}
D_x f_i(\vec{v})\\
D_y f_i(\vec{v})\end{array} \right) = \sum_{T \in N_{\vec{v}}} \nabla_T f_i/|N_{\vec{v}}|
\end{equation}

\noindent where $N_{\vec{v}}$ is the set of all faces around the vertex $\vec{v}$. Hence, the Beltrami coefficient $\mu_f(\vec{v})$ can be computed by:
\begin{equation}
\mu_f(\vec{v}) = \frac{(D_x f_1(\vec{v}) - D_y f_2(\vec{v})) + i(D_x f_2(\vec{v}) + D_y f_1(\vec{v}))}{(D_x f_1(\vec{v}) + D_y f_2(\vec{v})) + i(D_x f_2(\vec{v}) - D_y f_1(\vec{v}))}
\end{equation}

\noindent{\it 2. Computation of the BHF Reconstruction}

For the BHF reconstruction algorithm, the most important step is the computation of the variation $V(f^{\mu},\nu)$ of $f^{\mu}$ under the variation of $\mu$. We will discuss the computation of $V(f^{\mu},\nu)$ for $D=\mathbb{D}$. The computation for $D=\mathbb{S}^2 \equiv \overline{\mathbb{C}}$ is similar. From Equation \ref{BHFiteration} and \ref{BHFiteration2},
\begin{equation}
\begin{split}
V({f}^{\mu},\nu)&(w) =\int_D K(z,w) \,dx\,dy, \nonumber
\end{split}
\end{equation}

\noindent where
\begin{equation}
\begin{split}
K(z,w) &=-\frac{f^\mu(w)(f^\mu(w)-1)}{\pi}\\ \nonumber
&\left(\frac{\nu(z)((f^\mu)_z(z))^2}{f^\mu(z)(f^\mu(z)-1)(f^\mu(z)-f^\mu(w))}+ \frac{\overline{\nu(z)}(\overline{(f^\mu)_z(z)})^2}{\overline{f^\mu(z)}(1-\overline{f^\mu(z)})(1-\overline{f^\mu(z)}f^\mu(w))}\right)
\end{split}
\end{equation}

Now, $f^{\mu}$ and $\nu$ are both defined on each vertex $\vec{v}$. Also, $(f^\mu)_z (\vec{v})$ can be approximated as:
\begin{equation}
(f^\mu)_z (\vec{v}) \approx \frac{(D_x f_1(\vec{v}) - D_y f_2(\vec{v})) + i(D_x f_2(\vec{v}) + D_y f_1(\vec{v}))}{2}.
\end{equation}

For each pair of vertices $(\vec{v},\vec{w})$, $K(\vec{v},\vec{w})$ can be easily approximated. In case $K(\vec{v},\vec{w})$ is singular, we set $K(\vec{v},\vec{w})=0$. Now, for each vertex $\vec{v}$, we define $A_{\vec{v}}$ as
\begin{equation}
A_{\vec{v}} = \sum_{T\in N_{\vec{v}}} Area(T)/N_T,
\end{equation}

\noindent where $N_T$ is the number of vertices on $T$. That is, $N_T=3$ if $T$ is a triangle and $N_T=4$ if $T$ is a rectangle. Then, $V(f^{\mu},\nu)$ can be approximated by:
\begin{equation}
V(f^{\mu},\nu)(\vec{w}) = \sum_{\vec{v}} K(\vec{v},\vec{w}) A_{\vec{v}}
\end{equation}

\begin{center}
\underline{$\mathbf{II.\ Proof\ of\ Theorem\ \ref{thm:BHFD}}:$}
\end{center}

To prove the theorem, we need the following lemma:

\bigskip

\begin{lem}\label{thm:lemextension}
Let $f\colon \D \to \D$ be a diffeomorphism of the unit disk fixing $0$ and $1$ and satisfying the Beltrami equation $f_{\overline{z}}=\mu f_z$ with $\mu$ defined on $\D$. Let $\tilde{f}$ be the extension of $f$ to $\overline{\C}$ defined as
\begin{equation}
\tilde{f}(z)=
\begin{cases}
f(z),  & \mbox{if }|z| \leq 1,\\
\frac{1}{\overline{f(1/\overline{z})}}, & \mbox{if }|z|>1.
\end{cases}
\end{equation}
Then $\tilde{f}$ satisfies the Beltrami equation
\begin{equation}\label{eqt:extended_Beltrami}
\tilde{f}_{\overline{z}}=\tilde{\mu} \tilde{f}_z
\end{equation}
on $\overline{\C}$, where the Beltrami coefficient $\tilde{\mu}$ is defined as
\begin{equation}\label{eqt:extended_Beltrami2}
\tilde{\mu}(z)=
\begin{cases}
\mu(z), & \mbox{if }|z| \leq 1,\\
\frac{z^2}{\overline{z}^2}\overline{\mu(1/\overline{z})}, & \mbox{if }|z|>1.
\end{cases}
\end{equation}
\end{lem}

\bigskip

\begin{proof}
We need to prove $\tilde{f}$ satisfies the Beltrami equation:
\begin{equation}
\tilde{f}_{\overline{z}}=\tilde{\mu} \tilde{f}_z
\end{equation}
Clearly, $\tilde{f}$ satisfies equation (\ref{eqt:extended_Beltrami}) on $\D$. Outside $\D$, we consider $f$ and $\tilde{f}$ as functions in $z$ and $\overline{z}$.
Note that:
\begin{equation}
\frac{\partial}{\partial z}\overline{f({z},\overline{z})} = \overline{\frac{\partial}{\partial \overline{z}}{f({z},\overline{z})}}
\end{equation}
We have:
\begin{equation}\label{eqt1st}
\begin{split}
\frac{\partial\tilde{f}(z,\overline{z})}{\partial z} &=\frac{\partial}{\partial z}\frac{1}{\overline{f(1/\overline{z},1/z)}} =-\overline{f(1/\overline{z},1/z)}^{-2} \frac{\partial}{\partial z}\overline{f(1/\overline{z},1/z)} \\
    &=-\overline{f(1/\overline{z},1/z)}^{-2}\overline{\frac{\partial}{\partial\overline{z}}{f(1/\overline{z},1/z)}}=-\overline{f(1/\overline{z},1/z)}^{-2} \overline{(-1/\overline{z}^{2}) f_z(1/\overline{z},1/z)} \\
    &=z^{-2} \overline{f(1/\overline{z},1/z)}^{-2} \overline{f_z(1/\overline{z},1/z)}.
\end{split}
\end{equation}

Also,
\begin{equation}
\begin{split}
\frac{\partial\tilde{f}(z,\overline{z})}{\partial \overline{z}}
    &=\frac{\partial}{\partial \overline{z}}\frac{1}{\overline{f(1/\overline{z},1/z)}} =-\overline{f(1/\overline{z},1/z)}^{-2}\frac{\partial}{\partial \overline{z}}\overline{f(1/\overline{z},1/z)} \\
    &=-\overline{f(1/\overline{z},1/z)}^{-2} \overline{\frac{\partial}{\partial z}f(1/\overline{z},1/z)} =-\overline{f(1/\overline{z},1/z)}^{-2} \overline{(-1/z^{2}) f_{\overline{z}}(1/\overline{z},1/z)}\\
    &=\overline{z}^{-2} \overline{f(1/\overline{z},1/z)}^{-2} \overline{f_{\overline{z}}(1/\overline{z},1/z)}=\overline{z}^{-2} \overline{f(1/\overline{z},1/z)}^{-2} \overline{\mu(1/\overline{z})} \overline{f_z(1/\overline{z},1/z)}
\end{split}
\end{equation}

Now from Equation \ref{eqt1st},
\begin{equation}
\begin{split}
\overline{f_z(1/\overline{z},1/z)}
    &=z^{2} \overline{f(1/\overline{z},1/z)}^{2}\frac{\partial\tilde{f}(z,\overline{z})}{\partial z} .
\end{split}
\end{equation}

Thus, we have,
\begin{equation}
\begin{split}
\frac{\partial\tilde{f}(z,\overline{z})}{\partial \overline{z}}
    &=\overline{z}^{-2} \overline{f(1/\overline{z},1/z)}^{-2} \overline{\mu(1/\overline{z})} \overline{f_z(1/\overline{z},1/z)}\\
     &=\overline{z}^{-2} \overline{f(1/\overline{z},1/z)}^{-2} \overline{\mu(1/\overline{z})} z^{2}\overline{f(1/\overline{z},1/z)}^{2}\frac{\partial\tilde{f}(z,\overline{z})}{\partial z} \\
    &=\frac{z^2}{\overline{z}^2} \overline{\mu(1/\overline{z})} \frac{\partial \tilde{f}(z,\overline{z})}{\partial z}=\tilde{\mu}(z) \frac{\partial\tilde{f}(z,\overline{z})}{\partial z}.
\end{split}
\end{equation}
\end{proof}

\paragraph{Proof of Theorem \ref{thm:BHFD}}
According to Quasiconformal Teichmuller Theory, there is a one-to-one correspondence between the set of quasiconformal homeomorphisms of $\overline{\C}$ fixing 3 points and the set of smooth complex-valued functions $\mu$ on $\overline{\D}$ for which $\sup |\mu|=k<1$. If a diffeomorphism $f$ on $\overline{\C}$ satisfies equations (\ref{eqt:extended_Beltrami})(\ref{eqt:extended_Beltrami2}), then $1/\overline{f(1/\overline{z})}$ also satisfies the same equation. By the uniqueness of the solution according to Theorem \ref{thm:BHFC}, we must have $f(z)=1/\overline{f(1/\overline{z})}$. On $\partial\D$, $z=1/\overline{z}$. This implies $f(z)=1/\overline{f(z)}$, and hence $|f(z)|=1$ on $\partial\D$. Therefore, by restricting the solution of equation (\ref{eqt:extended_Beltrami}) on $\overline{\C}$ fixing $0$, $1$ and $\infty$ to $\D$, we can get a diffeomorphism of $\D$ fixing $0$ and $1$. Equation (\ref{eqt:flow}) can then be applied to get the variational formula $V(f^{\mu},\nu)$ of $f^{\mu}$ under the variation $\nu$ of $\mu$. To get $V(f^{\mu},\nu)$, we evaluate the integral in equation (\ref{eqt:flow}).
\begin{equation}
\begin{split}
&\int_\C \frac{\nu(z)((f^\mu)_z(z))^2}{f^\mu(z)(f^\mu(z)-1)(f^\mu(z)-f^\mu(w))}\,dx\,dy \\
    =&\int_\D \frac{\nu(z)((f^\mu)_z(z))^2}{f^\mu(z)(f^\mu(z)-1)(f^\mu(z)-f^\mu(w))}\,dx\,dy
        +\int_{\C\backslash\D} \frac{\nu(z)((f^\mu)_z(z))^2}{f^\mu(z)(f^\mu(z)-1)(f^\mu(z)-f^\mu(w))}\,dx\,dy \\
\end{split}
\end{equation}
Now, outside the disk $\D$,
\begin{equation}
\nu(z)=\frac{z^2}{\overline{z}^2} {\overline{\nu(1/\overline{z})}} \ \ \ \ \ \mathrm{and} \ \ \ \ \ \frac{\partial f(z)}{\partial z} = z^{-2} \overline{f(1/\overline{z},1/z)}^{-2} \overline{f_z(1/\overline{z},1/z)}
\end{equation}

We have:
\begin{equation}
\begin{split}
&\int_\C \frac{\nu(z)((f^\mu)_z(z))^2}{f^\mu(z)(f^\mu(z)-1)(f^\mu(z)-f^\mu(w))}\,dx\,dy \\
    =&\int_\D \frac{\nu(z)((f^\mu)_z(z))^2}{f^\mu(z)(f^\mu(z)-1)(f^\mu(z)-f^\mu(w))}\,dx\,dy
        +\int_{\C\backslash\D} \frac{\nu(z)((f^\mu)_z(z))^2}{f^\mu(z)(f^\mu(z)-1)(f^\mu(z)-f^\mu(w))}\,dx\,dy \\
    =&\int_\D \frac{\nu(z)((f^\mu)_z(z))^2}{f^\mu(z)(f^\mu(z)-1)(f^\mu(z)-f^\mu(w))}\,dx\,dy
        +\int_{\C\backslash\D} \frac{(z^2/\overline{z}^2) \overline{\nu(1/\overline{z})}((f^\mu)_z(z))^2}{f^\mu(z)(f^\mu(z)-1)(f^\mu(z)-f^\mu(w))}\,dx\,dy \\
    =&\int_\D \frac{\nu(z)((f^\mu)_z(z))^2}{f^\mu(z)(f^\mu(z)-1)(f^\mu(z)-f^\mu(w))}\,dx\,dy\\
       & +\int_\D \frac{(z^2/\overline{z}^2) \overline{\nu(z)}((f^\mu)_z(1/\overline{z}))^2} {\overline{f^\mu(1/\overline{z})}^{-1}(\overline{f^\mu(1/\overline{z})}^{-1}-1)(\overline{f^\mu(1/\overline{z})}^{-1}-f^\mu(w))}\frac{1}{|z|^4}\,dx\,dy\\
    =& \int_\D \frac{\nu(z)((f^\mu)_z(z))^2}{f^\mu(z)(f^\mu(z)-1)(f^\mu(z)-f^\mu(w))}\,dx\,dy +\int_\D \frac{ \overline{\nu(z)}(\overline{(f^\mu)_z(z)})^2} {\overline{f^\mu(z)}(1-\overline{f^\mu(z)})(1-\overline{f^\mu(z)}f^\mu(w))}\,dx\,dy
\end{split}
\end{equation}

Substituting Equation 22 into Equation \ref{eqt:flow}, we get an integral flow equation on $\D$, which is given by
\begin{equation}
\begin{split}
&V(f^{\mu}, \nu)(w)=-\frac{f^\mu(w)(f^\mu(w)-1)}{\pi} \\
&\left(\int_\D \frac{\nu(z)((f^\mu)_z(z))^2}{f^\mu(z)(f^\mu(z)-1)(f^\mu(z)-f^\mu(w))}\,dx\,dy +\int_\D \frac{ \overline{\nu(z)}(\overline{(f^\mu)_z(z)})^2} {\overline{f^\mu(z)}(\overline{1-f^\mu(z)})(1-\overline{f^\mu(z)}f^\mu(w))}\,dx\,dy.\right).
\end{split}
\end{equation}
\hspace*{\fill}~\QED\par\endtrivlist\unskip

{
\bibliographystyle{ieee}
\bibliography{egbib}
}
\end{document}